\documentclass[onefignum,onetabnum]{siamonline220329}

\usepackage{amsfonts}
\usepackage{amsmath}
\usepackage{graphicx}
\usepackage{float}
\usepackage{amssymb}
\usepackage{amsfonts}
\usepackage{xcolor}
\usepackage{amsmath}
\usepackage{bm}
\usepackage{cleveref}
\usepackage{diagbox}
\usepackage{stmaryrd}
\usepackage{mathtools}
\usepackage[normalem]{ulem}
\usepackage{overpic}
\usepackage{algorithm,algpseudocode}

\Crefname{subsection}{Section}{Sections}

\DeclareMathOperator{\tr}{tr}
\DeclareMathOperator{\Var}{Var}
\DeclareMathOperator{\rank}{rank}
\DeclareMathOperator{\range}{range}
\DeclareMathOperator{\vect}{vec}

\renewcommand*{\eqref}[1]{(\ref{#1})}
\newcommand{\tK}{\widetilde{\bm K}}
\newcommand{\bOmega}{\bm{\Omega}}
\newcommand{\R}{\mathbb{R}}
\newcommand{\Frob}{\textup{F}}
\newcommand{\Tr}{\mathrm{Tr}}
\newcommand{\HS}{\mathrm{HS}}
\newcommand{\op}{\mathrm{op}}
\renewcommand{\d}{\,\mathrm{d}}
\newcommand{\N}{\mathbb{N}}

\newcommand{\rev}[1]{\textcolor{black}{#1}}
\usepackage{stmaryrd}

\usepackage{enumitem}
\setlist[enumerate]{leftmargin=.5in}
\setlist[itemize]{leftmargin=.5in}

\newsiamremark{remark}{Remark}
\newsiamremark{hypothesis}{Hypothesis}
\crefname{hypothesis}{Hypothesis}{Hypotheses}
\newsiamthm{claim}{Claim}

\headers{Randomized Nystr{\"o}m approximation of non-negative self-adjoint operators}{D.~Persson, N.~Boull\'e, and D.~Kressner}

\title{Randomized Nystr{\"o}m approximation of non-negative self-adjoint  operators\thanks{Submitted to the editors \today{}.
        \funding{David Persson and Daniel Kressner were supported by the SNSF research project \textit{Fast algorithms from low-rank updates}, grant number: 200020\_178806. Nicolas Boull\'e was supported by an INI-Simons Postdoctoral Research Fellowship and the Office of Naval Research (ONR), under grant N00014-23-1-2729.}}}

\author{David Persson\thanks{Courant Institute of Mathematical Sciences, New York University, New York, NY 10012, USA \& Center for Computational Mathematics, Flatiron Institute, New York, NY 10010, USA (\email{dup210.nyu@edu}).}
    \and Nicolas Boull\'e\thanks{Department of Mathematics, Imperial College London, London SW7 2AZ, United Kingdom (\email{n.boulle@imperial.ac.uk}).}
    \and Daniel Kressner\thanks{Institute of Mathematics, EPF Lausanne, 1015 Lausanne, Switzerland. (\email{daniel.kressner@epfl.ch}).}}

\usepackage{amsopn}
\DeclareMathOperator{\diag}{diag}

\makeatletter
\newcommand{\proofpart}[2]{%
  \par
  \addvspace{\medskipamount}%
  \noindent\emph{#1. #2}\par\nobreak
  \addvspace{\smallskipamount}%
  \@afterheading
}
\makeatother

\begin{document}

\maketitle

\begin{abstract}
    The randomized singular value decomposition (SVD) has become a popular approach to computing cheap, yet accurate, low-rank approximations to matrices due to its efficiency and strong theoretical guarantees. Recent work by Boullé and Townsend (FoCM, 2023) presents an infinite-dimensional analog of the randomized SVD to approximate Hilbert-Schmidt operators. However, many applications involve computing low-rank approximations to symmetric positive semi-definite matrices. In this setting, it is well-established that the randomized Nystr{\"o}m approximation is usually preferred over the randomized SVD. This paper explores an infinite-dimensional analog of the Nystr{\"o}m approximation to compute low-rank approximations to non-negative self-adjoint trace-class operators. We present an analysis of the method and, along the way, improve the existing infinite-dimensional bounds for the randomized SVD. Our analysis yields bounds on the expected value and tail bounds for the Nystr{\"o}m approximation error in the operator, trace, and Hilbert-Schmidt norms. \rev{Numerical experiments on integral operators arising from Gaussian process sampling and Bayesian inverse problems are used to validate the proposed infinite-dimensional Nystr\"om algorithm.}
\end{abstract}



\section{Introduction}

Randomized techniques~\cite{gittensmahoney,rsvd,nakatsukasa2020fast,tropp2017fixed} are becoming increasingly popular for computing low-rank approximations to matrices.
Most notably, the \emph{randomized singular value decomposition} (SVD)~\cite{rsvd} has evolved into one of the primary choices, due to its simplicity, performance, and reliability. In its most basic form, the randomized SVD performs the approximation $\bm{Q} \bm{Q}^* \bm{A}\approx \bm{A}$, where $\bm{Q}$ is an orthonormal basis for the range of $\bm{A} \bm{\Omega}$, with $\bm{\Omega}$ being a tall and skinny random sketch matrix. In many applications of low-rank approximation,  such as $k$-means clustering \cite{wang2019scalable}, PCA \cite{zhang2008improved}, and Gaussian process regression \cite{gardner2018gpytorch}, it is known that $\bm{A}$ is symmetric positive semi-definite (SPSD). In this case, one usually prefers the so-called \emph{randomized Nystr{\"o}m approximation} \cite{gittensmahoney}
\begin{equation}\label{eq:nystrom}
    \widehat{\bm{A}} := \bm{A} \bm{\Omega}(\bm{\Omega}^* \bm{A} \bm{\Omega})^{\dagger} \bm{\Omega}^*\bm{A} \approx \bm{A},
\end{equation}
where $\bm{\Omega}$ is, again, a random sketch matrix and $\dagger$ denotes the Moore-Penrose pseudoinverse of a matrix. This approximation has received significant attention in the literature \cite{gittensmahoney,persson2022randomized,tropp2017fixed} and, like the randomized SVD, it enjoys strong theoretical guarantees. \rev{With the same number of matrix-vector products, the randomized Nyström approximation is typically significantly more accurate than the randomized SVD when the matrix has fast decaying singular values. Additionally, the Nyström method requires only a single pass over the matrix, compared to two passes for the randomized SVD, enabling all matrix-vector products to be performed in parallel. }

Recently, Boull\'e and Townsend~\cite{boulle2021a,boulle2022learning} generalized the randomized SVD from matrices to Hilbert-Schmidt operators. Subsequent works~\cite{boulle2022parabolic,boulle2023elliptic} employed this infinite-dimensional generalization of the randomized SVD to learn Green's functions associated with an elliptic or parabolic partial differential equations (PDE) from a few solutions of the PDE. This approach uses hierarchical low-rank techniques and exploits the fact that Green's functions are smooth away from the diagonal and therefore admit accurate off-diagonal low-rank approximations~\cite{bebendorf2008hierarchical,bebendorf2003existence}. The connection between numerical linear algebra and operator learning~\cite{boulle2023mathematical} motivates the development of infinite-dimensional analog of low-rank approximation algorithms, such as the randomized Nystr{\"o}m approximation, to derive sample complexity estimates of neural operator methods~\cite{kovachki2023neural,li2020fourier,lu2021learning} in scientific machine learning. Other applications, like Gaussian process regression and Support Vector Machines~\cite{eriksson2018,gardner2018gpytorch,pang2019,vapnik,williamsbarber,williams2000using}, involve integral operators that feature positive and \emph{globally} smooth kernels. In turn, the operator is not only self-adjoint and positive but it also allows for directly applying low-rank approximation, without the need to resort to hierarchical techniques. Given existing results on matrices, it would be sensible to use an infinite-dimensional extension of the randomized Nystr{\"o}m approximation in such situations.

\begin{remark}
    The use of the term \emph{Nystr{\"o}m approximation} in some of the recent literature and this work deviates somewhat from its original meaning, where the matrix $\bm{\Omega}$ in \eqref{eq:nystrom} is required to consist of columns of the identity matrix $\bm{I}$ and, in turn, $\bm{A} \bm{\Omega}$ consists of selected columns of $\bm{A}$. This form of Nystr{\"o}m approximation frequently appears in applications involving kernel matrices; see, e.g.,~\cite{bach2002kernel,burt2019rates,derezinski2020improved,drineas2005nystrom,williams2000using}. Following~\cite{gittensmahoney}, we interpret the term Nystr{\"o}m approximation in a broader sense that allows for other choices of matrices $\bOmega$, including Gaussian random matrices.
\end{remark}

\subsection{Contributions}

In this work, we present and analyze an infinite-dimensional extension of the randomized Nystr{\"o}m approximation for computing low-rank approximations to self-adjoint, positive, trace class operators. A difficulty in the design and analysis of our extension is that existing results of the finite-dimensional case~\eqref{eq:nystrom} usually assume that the columns of $\bm{\Omega}$ are drawn from a \emph{standard} Gaussian multivariate distribution $\mathcal{N}(\bm{0},\bm{I})$, which does not have a practically meaningful infinite-dimensional analog. Following \cite{boulle2022learning}, the columns of $\bm{\Omega}$ are replaced with random fields drawn from a Gaussian process $\mathcal{GP}(0, K)$, which is an infinite-dimensional analog of a \emph{non-standard} multivariate Gaussian distribution.

Therefore, to analyze the infinite-dimensional generalization of the Nystr{\"o}m approximation, we will proceed through an analysis of the finite-dimensional Nystr{\"o}m approximation \eqref{eq:nystrom} when the columns of $\bm{\Omega}$ are drawn from a non-standard Gaussian distribution $\mathcal{N}(\bm{0},\bm{K})$, for some general symmetric positive semi-definite matrix $\bm{K}$. We derive expectation and probability bounds between the matrix and its approximant in the Frobenius, spectral, and nuclear norms in \cref{theorem:expectation_frob,theorem:expectation_nystrom,theorem:tail_nystrom}. Then, we use continuity arguments to obtain bounds on the infinite-dimensional generalization of~\eqref{eq:nystrom} in \cref{theorem:inf_probabilistic}. As a byproduct of our analysis, we improve the analysis of the infinite-dimensional analog of the randomized SVD. In particular, unlike the bounds presented in~\cite{boulle2022learning}, our improved bounds coincide with the bounds obtained by Halko, Martinsson, and Tropp~\cite{rsvd} in the finite-dimensional case when $\bm{K}$ is chosen as the identity matrix.

\rev{A key benefit of the proposed framework is that once a low-rank approximation of the operator is computed, one can obtain a low-rank approximation to \emph{any} discretization of the operator at low additional cost. Clearly, an actual numerical implementation of our Nystr\"om approximation also requires discretization to realize the action of the operator $\mathcal{A}$ on a function. However, this action can be approximated adaptively and to high accuracy, say up to machine precision, so that the error is negligible.}

\subsection{Organization of the paper}
The paper is organized as follows. In \Cref{section:finitedimensions} we present the setting for the Nystr{\"o}m approximation in finite dimensions, give an overview of Gaussian random vectors, and provide the analysis of the Nystr{\"o}m approximation for matrices with correlated Gaussian sketches. We then present the generalization of the Nystr{\"o}m approximation to Hilbert-Schmidt operators in \Cref{section:infinitedimensions}. Next, we conclude with numerical experiments in \Cref{section:experiments}. \rev{Our analysis relies on properties of random matrices with Gaussian random entries. The statements and proofs of these properties are deferred to \Cref{sec_prop_Gaussian}.}

\section{Nyström approximation in finite dimensions with correlated Gaussian sketches}\label{section:finitedimensions}

In this section, we present an analysis of the randomized Nyström approximation defined in \eqref{eq:nystrom}, where the columns of $\bm{\Omega}$ are independent and identically distributed (i.i.d.) random vectors with distribution $\mathcal{N}(\bm{0},\bm{K})$ for some SPSD covariance matrix $\bm K$. We derive bounds in the Frobenius, spectral, and nuclear norms, which we denote respectively by
$\|\bm A\|_\Frob=(\sum_{i=1}^n \sigma_i^2)^{1/2}$, $\|\bm A\|_2=\sigma_1$, and $\|\bm A\|_*= \sum_{i=1}^n \sigma_i$, where $\sigma_1\geq \sigma_2 \geq\ldots\geq\sigma_n$ are the singular values of $\bm A$. More generally, $\|\bm A\|_{(s)}=(\sum_{i=1}^n\sigma_i^s)^{1/s}$ denotes the Schatten $s$-norm, \rev{which coincides with the Frobenius norm when $s=2$, the nuclear norm when $s=1$, and the spectral norm when $s = \infty$.}

We consider a symmetric positive semi-definite (SPSD) matrix $\bm{A} \in \mathbb{R}^{n \times n}$ with rank greater than $k\geq 1$ and eigenvalue decomposition $\bm A=\bm U \bm \Sigma \bm U^*$, partitioned as
\begin{equation} \label{eq:decompA}
    \bm{A} = \begin{bmatrix} \bm{U}_1 & \bm{U}_2 \end{bmatrix} \begin{bmatrix} \bm{\Sigma}_1 & \bm{0} \\ \bm{0} & \bm{\Sigma}_2 \end{bmatrix} \begin{bmatrix} \bm{U}_1^* \\ \bm{U}_2^* \end{bmatrix}, \quad \bm{\Sigma}_1 = \diag(\sigma_1,\ldots,\sigma_k), \quad \bm{\Sigma}_2 = \diag(\sigma_{k+1},\ldots,\sigma_n),
\end{equation}
where $\bm{U}_1 \in \mathbb{R}^{n \times k}$ contains the first $k$ columns of $\bm U$ and $\sigma_1 \geq \cdots \geq \sigma_n \geq 0$ are the eigenvalues of $\bm{A}$. For an oversampling parameter $p$ and an arbitrary sketch matrix $\bm \Omega \in \R^{n\times (k+p)}$, we denote the Nystr{\"o}m approximation of $\bm A$ defined in~\eqref{eq:nystrom} by $\widehat{\bm{A}}$.

\subsection{Distribution of Gaussian sketches} \label{sec_dist_Gauss}

To introduce some notation and basic background material, we let \[
    \bOmega_1 = \bm{U}_1^* \bOmega \in \R^{k\times(k+p)}, \quad
    \bOmega_2 = \bm{U}_2^* \bOmega\in \R^{(n-k)\times (k+p)},                                                                                                                                                                             \]
and, to analyze their distribution, we partition the matrix $\tK = \bm{U}^* \bm{K} \bm{U}$ as
\begin{equation} \label{eq:deftildek}
    \tK = \begin{bmatrix}
        \bm{U}_1^* \bm{K} \bm{U}_1 & \bm{U}_1^* \bm{K} \bm{U}_2 \\
        \bm{U}_2^* \bm{K} \bm{U}_1 & \bm{U}_2^* \bm{K} \bm{U}_2
    \end{bmatrix}
    =:
    \begin{bmatrix}
        \tK_{11} & \tK_{21}^* \\
        \tK_{21} & \tK_{22}
    \end{bmatrix}, \quad \tK_{11} \in \R^{k\times k}.
\end{equation}
We assume that $\tK_{11}$ is invertible, which allows us to define the Schur complement $\tK_{22.1} = \tK_{22}-\tK_{21}\tK_{11}^{-1}\tK_{21}^*$. By well-known properties of Gaussian random vectors~\cite[Thm.~5.2]{multstat}, $\bm \omega\sim \mathcal{N}(0,\bm K)$ implies $\bm U^*\bm{\omega} \sim \mathcal{N}(0,\tK)$ and, thus, the marginal distributions
\[
    \bm{\omega}_1:=\bm U_1^*\bm{\omega} \sim \mathcal{N}(0,\tK_{11}), \quad
    \bm{\omega}_2=\bm U_2^*\bm{\omega} \sim \mathcal{N}(0,\tK_{22}).
\]
In particular, the columns of $\bm \Omega_1$ are i.i.d. as $N(\bm{0}, \tK_{11})$ or, in other words, $\bOmega_1 \rev{\sim} \smash{\tK_{11}^{1/2} \bm{X}}$ with the standard Gaussian matrix $\bm{X}\in \R^{k\times (k+p)}$. Because $\tK_{11}$ is invertible, this shows that, with probability one, $\bOmega_1$ has full rank and possesses a right inverse $\bOmega_1^\dagger$, which satisfies $\bm \Omega_1\bm\Omega_1^\dagger=\bm I_k$.
Finally, the \emph{conditional} probability distribution of $\bm{\omega}_2$ given $\bm{\omega}_1=\bm{x}_1$ is normal with mean $\tK_{21} \tK_{11}^{-1} \bm{x}_1$ and covariance matrix $\tK_{22.1}$~\cite[Thm.~5.3]{multstat}.

\subsection{Nyström method with correlated Gaussian sketches}

Similarly to the error analysis for the randomized SVD with correlated input vectors~\cite{boulle2021a,boulle2022learning}, our error bounds for the Nystr{\"o}m approximation depend on prior information contained in the representation~\eqref{eq:deftildek} of $\bm{K}$ with respect to the eigenvectors of $\bm A$.  This is measured through the following two quantities:
\begin{equation} \label{eq_def_cov_factor}
    \beta_k^{(\xi)} = \frac{\|\bm{\Sigma}_2^{1/2} \tK_{22.1}\bm{\Sigma}_2^{1/2}\|_{\xi}}{\|\bm{\Sigma}_2\|_{\xi}}\|\tK_{11}^{-1}\|_2, \quad \rev{\delta_k^{(\xi)} = \frac{\|\bm{\Sigma}_2^{1/2} \tK_{21}\tK_{11}^{-2}\tK_{21}^*\bm{\Sigma}_2^{1/2}\|_{\xi}}{\|\bm{\Sigma}_2\|_{\xi}},}
\end{equation}
where $\xi\in\{\Frob,2,*\}$ such that $\|\cdot\|_{\xi}$ dictates the choice of norm. 

\rev{Note that $\widetilde{\bm{K}}_{22.1}$ is the covariance matrix of $\bm{\Omega}_2$ conditioned on $\bm{\Omega}_1$. Therefore, $\beta_k^{(\xi)}$ is small if the random vectors in $\bm{\Omega}$ have a small variance in the direction of the bottom eigenvectors of $\bm{A}$, conditioned on the randomness in the direction of the top $k$ eigenvectors of $\bm{A}$. In other words, $\beta_k^{(\xi)}$ measures the randomness in the eigenvectors corresponding to the lower eigenvalues of $\bm{A}$, weighted by the bottom $n-k$ eigenvalues of $\bm{A}$.
Furthermore, as discussed in \cite{di2024general}, note that $\widetilde{\bm{K}}_{21} \widetilde{\bm{K}}_{11}^{-1} = \tan(\bm{U}_1,\bm{K}\bm{U}_1)$ is a matrix whose singular values are the tangents of the principal angles between $\range(\bm{U}_1)$ and $\range(\bm{K} \bm{U}_1)$ \cite[Theorem 3.1]{zhu2013angles}. Therefore, to make $\delta_k^{(\xi)}$ small we should choose $\bm{K}$ so that $\bm{U}_1$ is an approximate invariance subspace for $\bm{K}$. Hence, $\delta_k^{(\xi)}$ is a measure of how well $\range(\bm{K} \bm{U}_1)$ and $\range(\bm{U}_1)$ align, weighted with the bottom $n-k$ eigenvalues of $\bm{A}$. }

\rev{The rest of the section is devoted to proving error bounds for the Nyström approximation. The proofs of our error bounds are based on the existing structural bound
\begin{equation} \label{eq:structurebound}
    \|\bm{A} - \widehat{\bm{A}}\| \leq \|\bm{\Sigma}_2\| + \|(\bm{\Sigma}_2^{1/2}\bm{\bOmega}_2 \bOmega_1^{\dagger})^* \bm{\Sigma}_2^{1/2}\bOmega_2 \bOmega_1^{\dagger}\|,
\end{equation}
where $\bm\Omega_1^\dagger=\bOmega_1^*(\bOmega_1\bOmega_1^*)^{-1}$ is the right inverse of $\bm\Omega_1$, assuming that this matrix has full rank, and $\|\cdot\|$ denotes any unitarily invariant norm. This bound is classical, and when $\|\cdot\|$ is a Schatten norm this bound follows from \cite[Proposition 8.5 and p.38]{tropp2023randomized} or \cite[Theorem 9.1]{rsvd}. When $\|\cdot\|$ is any unitarily invariant norm this bound follows from~\cite[Lem.~3.13]{persson2022randomized} with $q = 1$ and the identity function $f:x\mapsto x$. \rev{Similar bounds for the Frobenius, nuclear, and spectral norms can be found in \cite{gittensmahoney}.}
Obtaining probabilistic bounds from~\eqref{eq:structurebound} requires the analysis of the term $\|(\bm{\Sigma}_2^{1/2}\bm{\bOmega}_2 \bOmega_1^{\dagger})^* \bm{\Sigma}_2^{1/2}\bOmega_2 \bOmega_1^{\dagger}\|$.} 

\subsection{Analysis in the nuclear and operator norms}\label{section:othernorms}
\rev{The following theorem states our main result concerning the approximation error of the Nystr{\"o}m approximation in the nuclear and operator norms. Similar results are derived in the Frobenius norm but are deferred to \Cref{section:frobeniusnorm} for clarity of exposition.}

\subsubsection{Expectation bounds}
\rev{We begin with with the statement and proof of an expectation bounds for the error $\|\bm{A} - \widehat{\bm{A}}\|$ in the operator and nuclear norms. }
\begin{theorem}[Expectation bound in spectral and nuclear norms]\label{theorem:expectation_nystrom}
     \rev{Let $\bm A$ be an $n\times n$ SPSD matrix, $2\leq k\leq \rank(\bm A)$ be a target rank, and $p\ge 2$ an oversampling parameter. Let $\bm \Omega$ be a random sketch matrix with the columns i.i.d. $\mathcal{N}(\bm{0},\bm{K})$, where the covariance matrix $\bm K$ is such that the matrix $\tK_{11}$ defined in \eqref{eq:deftildek} is invertible. Then, the Nystr{\"o}m approximation $\widehat{\bm{A}}$ satisfies
    \begin{subequations} \label{eq_exp_spect_nucl}
        \begin{align}
            \mathbb{E}[\|\bm{A}-\widehat{\bm{A}}\|_2] & \leq \left(1 + \frac{3k}{p-1} \beta_k^{(2)} + 3\delta_k^{(2)}\right)\|\bm{\Sigma}_2\|_2 + \frac{3e^2(k+p)}{p^2-1}\beta_k^{(*)}\|\bm{\Sigma}_2\|_*, \label{eq_exp_spect} \\
            \mathbb{E}[\|\bm{A}-\widehat{\bm{A}}\|_*] & \leq \left(1 + \frac{k}{p-1} \beta_k^{(*)} +\delta_k^{(*)}\right)\|\bm{\Sigma}_2\|_*. \label{eq_exp_nucl}
        \end{align}
    \end{subequations}}
\end{theorem}

\begin{proof}
    \rev{To prove the inequalities \cref{eq_exp_spect} and  \cref{eq_exp_nucl}, we use the $L_q$ norm of a random variable $Z$ defined as $\mathbb E^q(Z) = \mathbb E[|Z|^q]^{1/q}$. After taking expectation with respect to $\bm\Omega$ and applying H\"older's inequality, the second term in the bound~\eqref{eq:structurebound} for the operator norm satisfies
    \[\mathbb{E}[\|\bm{A}-\widehat{\bm{A}}\|_2] \leq \|\bm\Sigma_2\|_2+\mathbb{E}[\|\bm{\Sigma}_2^{1/2}\bOmega_2\bOmega_1^{\dagger}\|_2^2] = \|\bm\Sigma_2\|_2+\mathbb{E}^2[\|\bm{\Sigma}_2^{1/2}\bOmega_2\bOmega_1^{\dagger}\|_2]^2.\]
    To proceed from here, we recall that \Cref{sec_dist_Gauss} provides the conditional distribution of $\bOmega_2|\bOmega_1$ as $\bOmega_2|\bOmega_1 \sim \tK_{21} \tK_{11}^{-1}\bOmega_1+\tK_{22.1}^{1/2}\bm\Psi$, where $\bm\Psi$ is a standard Gaussian matrix. We use the conditional distribution of $\bOmega_2|\bOmega_1$, the triangle inequality for the $L_2$ norm twice, and $\bm \Omega_1\bm\Omega_1^\dagger=\bm I_k$ (which holds with probability $1$), to get
    \begin{align*}
        \mathbb{E}^2\big[\|\bm{\Sigma}_2^{1/2}\bOmega_2 \bOmega_1^{\dagger}\|_{2}\big]
         & = \mathbb{E}^2_{\bm\Omega_1}\big[\mathbb{E}^2_{\bm\Psi} \big[\|\bm{\Sigma}_2^{1/2}\tK_{21} \tK_{11}^{-1} \bOmega_1 \bOmega_1^{\dagger} + \bm{\Sigma}_2^{1/2} \tK_{22.1}^{1/2} \bm{\Psi}\bOmega_1^{\dagger}\|_{2}\mid \bm \Omega_1\big]\big] \\
         & \le \|\bm{\Sigma}_2^{1/2}\tK_{21} \tK_{11}^{-1}\|_{2} + \big(
        \mathbb{E}_{\bm\Omega_1} \big[\mathbb{E}_{\bm\Psi}\big[\| \bm{\Sigma}_2^{1/2} \tK_{22.1}^{1/2} \bm{\Psi}\bOmega_1^{\dagger}\|_{2}^2\mid\bOmega_1\big]\big] \big)^{1/2}.
    \end{align*}
    Then \cref{eq_exp_shift_two} leads to
    \[\mathbb{E}_{\bm\Psi}\big[\| \bm{\Sigma}_2^{1/2} \tK_{22.1}^{1/2} \bm{\Psi}\bOmega_1^{\dagger}\|_{2}^2\mid\bOmega_1\big]\leq \big(\|\bm{\Sigma}_2^{1/2} \tK_{22.1}^{1/2}\|_\Frob\|\bOmega_1^\dagger\|_2+\|\bm{\Sigma}_2^{1/2} \tK_{22.1}^{1/2}\|_2\|\bOmega_1^\dagger\|_\Frob\big)^2.\]
    After taking expectation with respect to $\bOmega_1$ and using Hölder's inequality we have
    \begin{align*}
        \big(
        \mathbb{E}_{\bm\Omega_1} \big[\mathbb{E}_{\bm\Psi}\big[\| \bm{\Sigma}_2^{1/2} \tK_{22.1}^{1/2} \bm{\Psi}\bOmega_1^{\dagger}\|_{2}^2\mid\bOmega_1\big]\big] \big)^{1/2}
         & \leq \|\bm{\Sigma}_2^{1/2} \tK_{22.1}^{1/2}\|_\Frob(\mathbb{E}_{\bOmega_1}\big[\|\bOmega_1^\dagger\|_2^2\big])^{1/2}     \\
         & \quad + \|\bm{\Sigma}_2^{1/2} \tK_{22.1}^{1/2}\|_2(\mathbb{E}_{\bOmega_1}\big[\|\bOmega_1^\dagger\|_\Frob^2\big])^{1/2}.
    \end{align*}
    We then apply \cref{lemma:inversewishart_norm_expectation} to obtain the following inequality
    \begin{align*}
        \mathbb{E}^2[\|\bm{\Sigma}_2^{1/2}\bOmega_2 \bOmega_1^{\dagger}\|_{2}] & \leq
        \frac{\|\tK_{11}^{-1/2}\|_\Frob}{\sqrt{p-1}}\|\bm{\Sigma}_2^{1/2} \tK_{22.1}^{1/2}\|_2
        + \|\bm{\Sigma}_2^{1/2}\tK_{21} \tK_{11}^{-1}\|_{2}                                                                                                                                                                                               \\
                                                                               & \quad + \frac{e\sqrt{k+p}}{\sqrt{p^2-1}}\|\tK_{11}^{-1}\|_2^{1/2}\|\bm{\Sigma}_2^{1/2} \tK_{22.1}^{1/2}\|_\Frob                                                          \\
                                                                               & \leq \sqrt{\frac{k}{p-1}\beta_k^{(2)}\|\bm{\Sigma}_2\|_2} + \sqrt{\delta_k^{(2)} \|\bm{\Sigma}_2\|_2} + \sqrt{\frac{e^2(k+p)}{p^2-1} \beta_k^{(*)} \|\bm{\Sigma}_2\|_*},
    \end{align*}}
    \rev{where we used $\|\widetilde{\bm{K}}_{11}^{-1/2}\|_\Frob \leq \sqrt{k}\|\widetilde{\bm{K}}_{11}^{-1}\|_2^{1/2}$ and the definitions of $\beta_k^{(2)}$, $\beta_k^{(*)}$, and $\delta_k^{(2)}$. We conclude the proof of \cref{eq_exp_spect} using the inequality $(a + b + c)^2 \leq 3(a^2 + b^2 + c^2)$.}

    \rev{The bound for the nuclear norm follows through a similar argument from the structural bound~\eqref{eq:structurebound} in the nuclear norm:
    \begin{align*}
        \mathbb{E}[\|\bm{A}-\widehat{\bm{A}}\|_*] & \leq \|\bm\Sigma_2\|_*+\mathbb{E}[\|\bm{\Sigma}_2^{1/2}\bOmega_2\bOmega_1^{\dagger}\|_\Frob^2]                                                                \\
                                                  &\rev{=} \|\bm\Sigma_2\|_*+\|\bm{\Sigma}_2^{1/2}\tK_{21} \tK_{11}^{-1}\|_\Frob^2 + \|\bm{\Sigma}_2^{1/2}\tK^{1/2}_{22.1}\|_\Frob^2 \frac{\tr(\tK_{11}^{-1})}{p-1} \\
                                                  & \leq \left(1+\frac{k}{p-1}\beta_k^{(*)} + \delta_k^{(*)}\right)\|\bm{\Sigma}_2\|_*,
    \end{align*}
    where we used \eqref{eq_exp_shift_Frob} for the equality and $\tr(\widetilde{\bm{K}}_{11}^{-1}) \leq k \|\widetilde{\bm{K}}_{11}^{-1}\|_2$ for the inequality. }
\end{proof}

\subsubsection{Tail bounds}
\rev{We now proceed with the deriving tail bounds for $\|\bm{A} - \widehat{\bm{A}}\|$ in the operator and nuclear norms. We begin with the proof of a concentration inequality on norms of shifted and rescaled Gaussian matrices.}

\begin{lemma} \label{lemma:tail bound_shifted_pnorm}
    \rev{Let $\bm{\Psi}$ be a standard Gaussian matrix and let $\bm{B}, \bm{C}$, $\bm{D}$ be fixed matrices of matching sizes. Let $s \geq 2$. Then
    \[
        \mathbb{P}\left\{\|\bm{B}+\bm{C}\bm{\Psi}\bm{D}\|_{(s)} \geq \mathbb{E}\left[\|\bm{B}+\bm{C}\bm{\Psi}\bm{D}\|_{(s)}\right] + \|\bm{C}\|_2 \|\bm{D}\|_2 u\right\} \leq e^{-u^2/2}.
    \]
    holds for every $u\geq 1$.}
\end{lemma}
\begin{proof}
    \rev{Given $h(\bm{X}) := \|\bm{B}+\bm{C}\bm{X}\bm{D}\|_{(s)}$, we have
    \[
        |h(\bm{X})-h(\bm{Y})| \leq \|\bm{C}(\bm{X}-\bm{Y})\bm{D}\|_{(s)} \leq \|\bm{C}\|_2 \|\bm{D}\|_2\|\bm{X}-\bm{Y}\|_{(s)} \leq \|\bm{C}\|_2 \|\bm{D}\|_2\|\bm{X}-\bm{Y}\|_{\Frob},
    \]
    where we used that the Frobenius norm is larger than any Schatten-$s$ norm for $s \geq 2$.
    In other words, $h$ is Lipschitz continuous with Lipschitz constant $\|\bm{C}\|_2 \|\bm{D}\|_2$. This allows us
    to apply a concentration result for functions of Gaussian matrices~\cite[Prop.~10.3]{rsvd}, which yields the statement of the lemma.}
\end{proof}

\rev{With \cref{lemma:tail bound_shifted_pnorm} at hand, we are ready to state and prove the tail bound. }
\begin{theorem}[Tail bound in spectral and nuclear norms]\label{theorem:tail_nystrom}
    \rev{Using the notation of \cref{theorem:expectation_frob} with an oversampling parameter $p\geq 4$, define
    \begin{equation} \label{eq_def_di}
        d_1 = \frac{3 k}{p+1}, \quad d_2 = e^2\frac{k+p}{(p+1)^2}.
    \end{equation}
    Then for any $u,t\geq 1$, each of the two inequalities
    \begin{subequations} \label{eq_tail_spect_nucl}
        \begin{align}
            \|\bm{A}-\widehat{\bm{A}}\|_2 & \leq \left(1+4\delta_k^{(2)}+4(d_1+d_2 u^2) t^2\beta_k^{(2)}\right)\|\bm{\Sigma}_2\|_2+4 d_2 t^2\beta_k^{(*)}\|\bm{\Sigma}_2\|_*, \label{eq_tail_spect} \\
            \|\bm{A}-\widehat{\bm{A}}\|_* & \leq \left(1 + 2 \delta_k^{(*)} + d_1 t^2\beta_k^{(*)}\right)\|\bm{\Sigma}_2\|_*+2 d_2 t^2 u^2 \beta_k^{(2)}\|\bm{\Sigma}_2\|_2, \label{eq_tail_nucl}
        \end{align}
    \end{subequations}
    holds with probability at least $1- 2t^{-p}-e^{u^2/2}$.}
\end{theorem}

\begin{proof}
    \rev{We begin by deriving the tail bound \cref{eq_tail_spect} in the spectral norm. To obtain~\cref{eq_tail_spect}, it suffices to derive a tail bound for the term $\|\bm \Sigma_2^{1/2}\bOmega_2\bOmega_1^\dagger\|_{(2)}^2$ in the structural bound~\eqref{eq:structurebound}. Using that $\bOmega_2|\bOmega_1 \sim \tK_{21} \tK_{11}^{-1}\bOmega_1+\tK_{22.1}^{1/2}\bm\Psi$ with a standard Gaussian matrix $\bm\Psi$, \cref{lemma:tail bound_shifted_pnorm} yields the following tail bound for $\|\bm \Sigma_2^{1/2}\bOmega_2\bOmega_1^\dagger\|_{2}$ conditioned on $\bOmega_1$:
    \begin{equation} \label{eq_tail_concentration_sp}
        \mathbb{P}\left\{\|\bm \Sigma_2^{1/2}\bOmega_2\bOmega_1^\dagger\|_{2} \geq \mathbb{E}\big[\|\bm \Sigma_2^{1/2}\bOmega_2\bOmega_1^\dagger\|_{2}\mid \bOmega_1\big] + \|\bm \Sigma_2^{1/2}\tK_{22.1}^{1/2}\|_2 \|\bOmega_1^\dagger\|_2 u\mid \bOmega_1\right\} \leq e^{-u^2/2}.
    \end{equation}
    Using \cref{lemma:shifted_4norm} and the triangle inequality, it holds that
    \begin{align*}
        \mathbb{E}\big[\|\bm \Sigma_2^{1/2}\bOmega_2\bOmega_1^\dagger\|_2\mid \bOmega_1\big]                                                                                                                                 \leq \|\bm{\Sigma}_2^{1/2}\tK_{21} \tK_{11}^{-1}\|_{2}+\|\bm{\Sigma}_2^{1/2}\tK_{22.1}^{1/2}\|_{\Frob} \|\bOmega_1^{\dagger}\|_{2}+ \|\bm{\Sigma}_2^{1/2}\tK_{22.1}^{1/2}\|_{2} \|\bOmega_1^{\dagger}\|_{\Frob}.
    \end{align*}
    We then consider the probability event $E_t$ that the Frobenius and spectral norms of $\bOmega_1^{\dagger}$ are well controlled:
    \[E_t=\left\{\|\bOmega_1^{\dagger}\|_\Frob \leq \sqrt{d_1\|\tK_{11}^{-1}\|_2} t, \quad
        \|\bOmega_1^{\dagger}\|_{2} \leq \sqrt{d_2 \|\tK_{11}^{-1}\|_2} t\right\},\]
    where \rev{$d_1$ and $d_2$ are defined in \eqref{eq_def_di} and } $\mathbb{P}(E_t^c)\leq 2 t^{-p}$ by \cref{lemma:omega1_tailbound}. Then, conditioned on $E_t$, by \cref{eq_tail_concentration_sp} and the inequality $\|\bm{x}\|_1\leq \sqrt{d}\|\bm{x}\|_2$ for $\bm{x}\in \R^d$ we have:
    \[\|\bm \Sigma_2^{1/2}\bOmega_2\bOmega_1^\dagger\|_{2}^2 \leq \left(4\delta_k^{(2)}+4(d_1+d_2 u^2) t^2\beta_k^{(2)}\right)\|\bm{\Sigma}_2\|_2+4 d_2 t^2\beta_k^{(*)}\|\bm{\Sigma}_2\|_*,\]
    with probability $\geq 1-e^{-u^2/2}$.} \rev{Similarly to the proof of~\cite[Thm.~10.8]{rsvd}, we remove the conditioning by a union bound using $\mathbb{P}(E_t^c)\leq 2 t^{-p}$ and conclude the proof of~\eqref{eq_tail_spect}. The proof of the tail bound for the nuclear norm~\cref{eq_tail_nucl} follows from a similar argument.}
\end{proof}

\begin{remark}[Connection with the randomized SVD]\label{remark:rsvd}
    \rev{Consider an arbitrary matrix $\bm{B} \in \mathbb{R}^{m \times n}$ with singular value decomposition
    \[
        \bm{B} = \begin{bmatrix} \bm{W}_1 & \bm{W}_2 \end{bmatrix} \begin{bmatrix} \bm{S}_1 & \bm{0}\\ \bm{0} & \bm{S}_2 \end{bmatrix} \begin{bmatrix} \bm{U}_1^* \\ \bm{U}_2^* \end{bmatrix},
    \]
    where $\bm{W}_1 \in \mathbb{R}^{m \times k}, \bm{U}_1 \in \mathbb{R}^{n \times k}$, and $\bm{S}_1 \in \mathbb{R}^{k \times k}$ is a diagonal matrix containing the $k$ largest singular values of $\bm{B}$.
    Letting $\bm{Q}$ be an orthonormal basis for $\range(\bm{B} \bm{\Omega})$, then $\bm{Q} \bm{Q}^* \bm{B}$ is the approximation attained by the (basic) randomized SVD. Setting $\bm{A} = \bm{B}^* \bm{B}$ and $\widehat{\bm{A}} = \bm{A} \bm{\Omega} (\bm{\Omega}^* \bm{A} \bm{\Omega})^{\dagger} \bm{\Omega}^* \bm{A}$ we have for any $s \geq 1$ that
    \begin{equation*}
        \|\bm{B} - \bm{Q} \bm{Q}^* \bm{B}\|_{(2s)}^2 = \|\bm{A} - \widehat{\bm{A}}\|_{(s)}.
    \end{equation*}
    Hence, obtaining a bound for the randomized SVD applied to $\bm{B}$ is equivalent to obtaining a bound for the Nystr{\"o}m approximation on $\bm{A}$; similar observations have been made in~\cite{gittensmahoney,tropp2023randomized}. Therefore, one can apply the results obtained earlier in this section to derive error bounds for the randomized SVD. As an example, \cref{theorem:expectation_nystrom} implies that
    \begin{equation}\label{eq:ourbound}
        \mathbb{E}\|\bm{B} - \bm{Q} \bm{Q}^* \bm{B}\|_\Frob^2 \leq \left(1 + \frac{k}{p-1}\tilde{\beta}_k +\tilde{\delta}_k\right)\|\bm{S}_2\|_\Frob^2,
    \end{equation}
    where $\tilde{\beta}_k = \tr(\bm{S}_2 \widetilde{\bm{K}}_{22.1} \bm{S}_2)\|\widetilde{\bm{K}}_{11}^{-1}\|_2/\|\bm{S}_2\|_\Frob^2$ and $\tilde{\delta}_k = \tr(\bm{S}_2 \widetilde{\bm{K}}_{21} \widetilde{\bm{K}}_{11}^{-2} \widetilde{\bm{K}}_{21}^*\bm{S}_2)/\|\bm{S}_2\|_\Frob^2$. 
    \rev{In comparison, by following the analysis of the prior work in \cite{boulle2021a} one can show that}
        \begin{equation}\label{eq:oldbound}
        \rev{\mathbb{E}\|\bm{B} - \bm{Q} \bm{Q}^* \bm{B}\|_\Frob^2 \leq \left(1 + \frac{k(k+p)}{p-1}\tilde{\gamma}_k\right)\|\bm{S}_2\|_\Frob^2,}
    \end{equation}}%
    \rev{where $\tilde{\gamma}_k = \tr(\bm{S}_2 \widetilde{\bm{K}}_{22} \bm{S}_2)\|\widetilde{\bm{K}}_{11}^{-1}\|_2/\|\bm{S}_2\|_\Frob^2.$ Since $\tilde{\beta}_k, \tilde{\delta}_k \leq \tilde{\gamma}_k$, the inequality \eqref{eq:ourbound} improves on \eqref{eq:oldbound}. Furthermore, \eqref{eq:ourbound} coincides with the standard randomized SVD bound~\cite[Thm.~10.5]{rsvd} when $\bm{K} = \bm{I}$, unlike \eqref{eq:oldbound}.}
\end{remark}

\subsection{Analysis in the Frobenius norm}\label{section:frobeniusnorm}
\rev{We now proceed with the error bounds for $\|\bm{A} -\widehat{\bm{A}}\|_\Frob$. 
These are technical, but important to derive error bounds in Hilbert-Schmidt norm in the continuous setting in \Cref{section:infinitedimensions}. The Hilbert-Schmidt norm controls the variance of the stochastic trace estimator of integral operators, see \cite{conthutch++}, and is equal to the $L^2$-norm error of the kernel of the integral operators. The key ideas of the analysis will be similar to the analysis in \Cref{section:othernorms}, but requires bounding $\|(\bm{\Sigma}_2^{1/2}\bm{\bOmega}_2 \bOmega_1^{\dagger})^* \bm{\Sigma}_2^{1/2}\bOmega_2 \bOmega_1^{\dagger}\|_\Frob = \|\bm{\Sigma}_2^{1/2}\bOmega_2 \bOmega_1^{\dagger}\|_{(4)}^2$, which is more complicated. }
\subsubsection{Expectation bound}
\rev{We begin with stating and proving the expectation bound. }
\begin{theorem}[Frobenius norm] \label{theorem:expectation_frob}
    \rev{Consider the setting of \cref{theorem:expectation_nystrom} with an oversampling parameter $p \geq 4$. 
    \rev{Define the following constants
    \begin{equation} \label{eq_const_c1_c2}
        c_1 = k \frac{(p-1)(k+1) + 2}{p(p-1)(p-3)},\quad
    c_2 = k \frac{k + p-1}{p(p-1)(p-3)},
    \end{equation}}
    Then, the Nystr{\"o}m approximation $\widehat{\bm{A}}$ satisfies
    \begin{equation} \label{eq_exp_Frob}
        \mathbb{E}[\|\bm{A}-\widehat{\bm{A}}\|_\Frob] \leq \left( 1+2 \delta_k^{(\Frob)} + 2\sqrt{c_1} \beta_k^{(\Frob)}  \right) \|\bm\Sigma_2\|_\Frob
        + 2\sqrt{c_2} \beta_k^{(*)}  \|\bm\Sigma_2\|_*.
    \end{equation}}
    
    
\end{theorem}


\begin{proof}
\rev{To prove inequality \cref{eq_exp_Frob} of~\cref{theorem:expectation_frob}, we once again use the $L_q$ norm of a random variable $Z$ defined as $\mathbb E^q(Z) = \mathbb E[|Z|^q]^{1/q}$. After taking expectation with respect to $\bm\Omega$ and applying H\"older's inequality, the second term in the bound~\eqref{eq:structurebound} for the Frobenius norm satisfies
\[ \mathbb{E}[\|(\bm{\Sigma}_2^{1/2}\bm{\bOmega}_2 \bOmega_1^{\dagger})^* \bm{\Sigma}_2^{1/2}\bm{\bOmega}_2 \bOmega_1^{\dagger}\|_\Frob]
    =  \mathbb{E}[\|\bm{\Sigma}_2^{1/2}\bOmega_2 \bOmega_1^{\dagger}\|_{(4)}^2]
    \leq \big( {\mathbb{E}^4[\|\bm{\Sigma}_2^{1/2}\bOmega_2 \bOmega_1^{\dagger}\|_{(4)}]}\big)^{2}.\]
To proceed from here, we, once again, recall that \Cref{sec_dist_Gauss} provides the conditional distribution of $\bOmega_2|\bOmega_1$ as $\bOmega_2|\bOmega_1 \sim \tK_{21} \tK_{11}^{-1}\bOmega_1+\tK_{22.1}^{1/2}\bm\Psi$, where $\bm\Psi$ is a standard Gaussian matrix.
Using the triangle inequality for the $L_4$ norm and $\bm \Omega_1\bm\Omega_1^\dagger=\bm I_k$ (which holds with probability $1$), we obtain that
\begin{equation} \label{eq:ref123}
    \begin{aligned}
        \mathbb{E}^4[\|\bm{\Sigma}_2^{1/2}\bOmega_2 \bOmega_1^{\dagger}\|_{(4)}]
         & = \mathbb{E}_{\bOmega_1,\bm{\Psi}}^4[\|\bm{\Sigma}_2^{1/2}\tK_{21} \tK_{11}^{-1}\bOmega_1 \bOmega_1^{\dagger} + \bm{\Sigma}_2^{1/2}\tK_{22.1}^{1/2}\bm\Psi \bOmega_1^{\dagger} \|_{(4)}] \\
         & \le \|\bm{\Sigma}_2^{1/2}\tK_{21} \tK_{11}^{-1}\|_{(4)} + \big(
        \mathbb{E}_{\bm\Omega_1} \big[\mathbb{E}_{\bm\Psi}\big[\| \bm{\Sigma}_2^{1/2} \tK_{22.1}^{1/2} \bm{\Psi}\bOmega_1^{\dagger}\|_{(4)}^4\mid\bOmega_1\big]\big] \big)^{1/4},
    \end{aligned}
\end{equation}
To bound the second term, we first apply~\cref{eq_exp_shift_four},
\[
    \mathbb{E}_{\bm\Psi}\big[\| \bm{\Sigma}_2^{1/2} \tK_{22.1}^{1/2} \bm{\Psi}\bOmega_1^{\dagger}\|_{(4)}^4\mid\bOmega_1\big] = \|\bm{\Sigma}_2^{1/2}\tK_{22.1}^{1/2}\|_{(4)}^4 \big(\|\bOmega_1^{\dagger}\|_{(4)}^4 + \|\bOmega_1^{\dagger}\|_{\Frob}^4\big) + \|\bm{\Sigma}_2^{1/2}\tK_{22.1}^{1/2}\|_{\Frob}^4 \|\bOmega_1^{\dagger}\|_{(4)}^4,
\]
and then take expectation with respect to $\bm\Omega_1$ using \rev{\cref{lemma:shifted_4norm} and inequalities in}~\cref{lemma:inversewishart_norm_expectation}:
\[
    \rev{\mathbb{E}_{\bOmega_1} \big[\mathbb{E}_{\bm\Psi}\big[\| \bm{\Sigma}_2^{1/2} \tK_{22.1}^{1/2} \bm{\Psi}\bOmega_1^{\dagger}\|_{(4)}^4\mid\bOmega_1\big]\big] \leq
    \left(c_1 \|\bm{\Sigma}_2^{1/2}\tK_{22.1}^{1/2}\|_{(4)}^4 + c_2  \|\bm{\Sigma}_2^{1/2}\tK_{22.1}^{1/2}\|_{\Frob}^4\right)\|\tK_{11}^{-1}\|_2^2,}
\]
where \rev{$c_1$ and $c_2$ are defined in \eqref{eq_const_c1_c2}.}
Inserting the inequality above into \eqref{eq:ref123} and using $\|\tK_{11}^{-1/2}\|_2  = \|\tK_{11}^{-1}\|_2^{1/2}$ gives
\begin{align*}
     & \mathbb{E}^4[\|\bm{\Sigma}_2^{1/2}\bOmega_2 \bOmega_1^{\dagger}\|_{(4)}] \leq \\
     & \quad\qquad \|\bm{\Sigma}_2^{1/2}\tK_{21} \tK_{11}^{-1}\|_{(4)} + \big( c_1
    \|\bm{\Sigma}_2^{1/2}\tK_{22.1}^{1/2}\|_{(4)}^4 + c_2  \|\bm{\Sigma}_2^{1/2}\tK_{22.1}^{1/2}\|_{\Frob}^4\big)^{1/4} \|\tK_{11}^{-1}\|_2^{1/2},
\end{align*}
Finally, inserting the covariance quality factors (see \cref{eq_def_cov_factor}) leads to the following inequality,
\[
    \big( \mathbb{E}^4[\|\bm{\Sigma}_2^{1/2}\bOmega_2 \bOmega_1^{\dagger}\|_{(4)}]\big)^2 \le
    \left( 2 \delta_k^{(\Frob)} + 2\sqrt{c_1} \beta_k^{(\Frob)}  \right) \|\bm\Sigma_2\|_{\Frob}
    + 2\sqrt{c_2} \beta_k^{(*)}  \|\bm\Sigma_2\|_*,
\]
where we used $(a + b)^2 \leq 2a^2 + 2b^2$ and the subadditivity of the square-root. This concludes the proof of \cref{eq_exp_Frob}.}
\end{proof}

\subsubsection{Tail bound} \label{sec_tail_proof}


\rev{We now proceed with stating and proving the tail bound for the Frobenius norm error. }

\begin{theorem}[Frobenius norm] \label{theorem:tail bound_frob}
    \rev{Consider the setting of \cref{theorem:expectation_nystrom} with an oversampling parameter $p \geq 4$.} \rev{Define $d_3 = \sqrt{k}d_2$, where $d_2$ is defined in \eqref{eq_def_di} and let $u,t\geq 1$. Then,}
    \rev{
    \begin{align}
        \|\bm{A}-\widehat{\bm{A}}\|_\Frob & \leq
        \|\bm{\Sigma}_2\|_\Frob + 4 \left(\delta_k^{(\Frob)} + t^2(d_1+d_3)\beta_k^{(\Frob)}\right)\|\bm\Sigma_2\|_\Frob + 4 t^2 d_3\beta_k^{(*)}\|\bm\Sigma_2\|_* \label{eq_tail_Frob} \\
                                          & \quad + 2 t^2 u^2 d_2\beta_k^{(2)}\|\bm\Sigma_2\|_2 \nonumber
    \end{align}
    holds with probability at least $1-3t^{-p}-e^{-u^2/2}$. }
\end{theorem}
\begin{proof}
\rev{To obtain~\cref{eq_tail_Frob}, by \eqref{eq:structurebound} it suffices to derive a tail bound for $\|\bm \Sigma_2^{1/2}\bOmega_2\bOmega_1^\dagger\|_{(4)}^2$.
Using that $\bOmega_2|\bOmega_1 \sim \tK_{21} \tK_{11}^{-1}\bOmega_1+\tK_{22.1}^{1/2}\bm\Psi$ with a standard Gaussian matrix $\bm\Psi$, \cref{lemma:tail bound_shifted_pnorm} yields the following tail bound for $\|\bm \Sigma_2^{1/2}\bOmega_2\bOmega_1^\dagger\|_{(4)}^2$ conditioned on $\bOmega_1$:
\begin{equation} \label{eq_tail_concentration}
    \mathbb{P}\left\{\|\bm \Sigma_2^{1/2}\bOmega_2\bOmega_1^\dagger\|_{(4)} \geq \mathbb{E}\big[\|\bm \Sigma_2^{1/2}\bOmega_2\bOmega_1^\dagger\|_{(4)}\mid \bOmega_1\big] + \|\bm \Sigma_2^{1/2}\tK_{22.1}^{1/2}\|_2 \|\bOmega_1^\dagger\|_2 u\mid \bOmega_1\right\} \leq e^{-u^2/2}.
\end{equation}
\rev{For a deterministic matrix $\bm{X}$ and random matrix $\bm{Y}$, combining Jensen's inequality and the triangle inequality yields $\mathbb{E}\|\bm{X} + \bm{Y}\|_{(4)} \leq \mathbb{E}^4\left[\|\bm{X} + \bm{Y}\|_{(4)}\right] \leq \|\bm{X}\|_{(4)} + \mathbb{E}^4\left[\|\bm{Y}\|_{(4)}\right]$. Hence, by~\cref{eq_exp_shift_four} we get}
\begin{gather*}
    \mathbb{E}\big[\|\bm \Sigma_2^{1/2}\bOmega_2\bOmega_1^\dagger\|_{(4)}\mid \bOmega_1\big]\leq \|\bm{\Sigma}_2^{1/2}\tK_{21} \tK_{11}^{-1}\|_{(4)} +\mathbb{E}_{\bm\Psi}^4\big[\|\bm\Sigma_2^{1/2}\tK_{22.1}^{1/2}\bm{\Psi}\bOmega_1^\dagger\|_{(4)}\mid\bOmega_1\big]\\
    \rev{=} \|\bm{\Sigma}_2^{1/2}\tK_{21} \tK_{11}^{-1}\|_{(4)}+\left(\|\bm{\Sigma}_2^{1/2}\tK_{22.1}^{1/2}\|_{(4)}^4 \big(\|\bOmega_1^{\dagger}\|_{(4)}^4 + \|\bOmega_1^{\dagger}\|_{\Frob}^4\big) + \|\bm{\Sigma}_2^{1/2}\tK_{22.1}^{1/2}\|_{\Frob}^4 \|\bOmega_1^{\dagger}\|_{(4)}^4\right)^{1/4}.
\end{gather*}
Using, again, the inequality $(a+b)^2\leq 2(a^2+b^2)$ and $\sqrt{a+b}\leq \sqrt{a}+\sqrt{b}$ and the subadditivity of the square-root, we have
\begin{align}
    \mathbb{E}\big[\|\bm \Sigma_2^{1/2}\bOmega_2\bOmega_1^\dagger\|_{(4)}\mid \bOmega_1\big]^2
     & \leq
    2\|\bm{\Sigma}_2^{1/2}\tK_{21} \tK_{11}^{-1}\|_{(4)}^2 + 2\|\bm{\Sigma}_2^{1/2}\tK_{22.1}^{1/2}\|_{(4)}^2 (\|\bOmega_1^{\dagger}\|_{(4)}^2 + \|\bOmega_1^{\dagger}\|_{\Frob}^2) \label{eq_bound_exp_omega_tail} \\
     & \quad + 2\|\bm{\Sigma}_2^{1/2}\tK_{22.1}^{1/2}\|_{\Frob}^2 \|\bOmega_1^{\dagger}\|_{(4)}^2. \nonumber
\end{align}
To control the norms of $\bOmega_1^\dagger$ in our bounds, we condition on the event $E_t$ that the following three inequalities are satisfied:
\[
    \|\bOmega_1^{\dagger}\|^2_\Frob \leq {d_1\|\tK_{11}^{-1}\|_2} t^2, \quad
    \|\bOmega_1^{\dagger}\|^2_{2} \leq {d_2 \|\tK_{11}^{-1}\|_2} t^2, \quad
    \|\bOmega_1^{\dagger}\|^2_{(4)} \leq {d_3 \|\tK_{11}^{-1}\|_2} t^2,\]
with the constants} \rev{$d_1, d_2,$ and $d_3$ defined in \eqref{eq_def_di} and in the statement of \Cref{theorem:tail bound_frob}. }
\cref{lemma:omega1_tailbound} \rev{together with the inequalities $\tr(\widetilde{\bm{K}}_{11}^{-1}) \leq k \|\widetilde{\bm{K}}_{11}^{-1}\|_2$ and $\|\widetilde{\bm{K}}_{11}^{-1}\|_\Frob \leq \sqrt{k} \|\widetilde{\bm{K}}_{11}^{-1}\|_2$} \rev{implies the following bound for $\mathbb{P}(E_t^c)$:
\begin{equation} \label{eq_event_et}
    \mathbb{P}(E_t^c)\leq 2t^{-(p+1)}+t^{-p}\leq 3 t^{-p}.
\end{equation}
Under the event $E_t$, \cref{eq_bound_exp_omega_tail} becomes
\begin{equation} \label{eq_exp_Et}
    \mathbb{E}\big[\|\bm \Sigma_2^{1/2}\bOmega_2\bOmega_1^\dagger\|_{(4)}\mid E_t\big]^2
    \leq
    2\left(\delta_k^{(\Frob)} + t^2(d_1+d_3)\beta_k^{(\Frob)}\right)\|\bm\Sigma_2\|_\Frob + 2 t^2 d_3\beta_k^{(*)}\|\bm\Sigma_2\|_*.
\end{equation}
Conditioning \cref{eq_tail_concentration} on $E_t$ and combining it with \cref{eq_exp_Et} yields
\[
    \|\bm \Sigma_2^{1/2}\bOmega_2\bOmega_1^\dagger\|_{(4)}^2 \leq 4\left(\delta_k^{(\Frob)} + t^2(d_1+d_3)\beta_k^{(\Frob)}\right)\|\bm\Sigma_2\|_\Frob + 4 t^2 d_3\beta_k^{(*)}\|\bm\Sigma_2\|_* + 2d_2\beta_k^{(2)}\|\bm\Sigma_2\|_2 u^2 t^2,
\]
with probability $\geq 1-e^{-u^2/2}$.} \rev{Once again, we remove the conditioning by a union bound using \cref{eq_event_et} and conclude the proof of~\cref{theorem:expectation_frob}.}
\end{proof}

\section{Nystr{\"o}m approximation in infinite dimensions}\label{section:infinitedimensions}
This section presents an infinite-dimensional extension of the randomized Nystr{\"o}m approximation. We begin by briefly introducing the concepts of quasimatrices, Hilbert-Schmidt operators, and Gaussian processes, which will be useful to generalize the bounds of~\Cref{section:finitedimensions} to operators between function spaces. While we consider operators between square-integrable functions as the most important application, the results of this section extend naturally to other separable Hilbert spaces following the theory in~\cite{functionaldatanalysis}.

\subsection{Quasimatrices}\label{section:quasimatrix}

For a bounded domain $D\subset \R^d$, $d\geq 1$, we consider the Hilbert space $L^2(D)$ of \rev{real-valued} square-integrable functions. Quasimatrices are a convenient way to represent and work with collections of functions or more, generally, elements of infinite-dimensional vector spaces; see, e.g.,~\cite{quasimatrix}. In particular, a \rev{linear operator} $Y:\R^m\to L^2(D)$ \rev{can be} expressed as the quasimatrix
\begin{equation}\label{eq:quasi}
    Y=
    \begin{bmatrix}
        y_1 & \cdots & y_m
    \end{bmatrix}, \quad y_i \in L^2(D)\rev{,}
\end{equation}
\rev{where $y_i = Y\bm{e}_i$ and $\bm{e}_i$ is the $i^{\text{th}}$ canonical vector in $\mathbb{R}^{m}$ for $i = 1,\ldots,m$. Hence, similar} to matrices, compositions of linear operators can be conveniently expressed by extending the usual matrix multiplication rules to quasimatrices. 
\rev{If $Y$ is the quasimatrix \eqref{eq:quasi}, its adjoint $Y^* :L^2(D) \to \R^m$ can also be viewed as a quasimatrix with the $m$ rows $\langle y_1, \cdot \rangle, \cdots, \langle y_m, \cdot \rangle : L^2(D) \to \R^m$, where $\langle\cdot,\cdot\rangle$ denotes the standard inner product in $L^2(D)$. Then, if $f \in L^2(D)$ we have
\begin{equation*}
    Y^* f = \begin{bmatrix} \langle y_1,f \rangle \\ \vdots\\ \langle y_m, f \rangle \end{bmatrix} \in \mathbb{R}^m.
\end{equation*}}
If $Z:\R^\ell\to L^2(D)$ is another \rev{quasimatrix}, then $Y^*Z$ \rev{yields a linear operator $\mathbb{R}^{\ell} \to \mathbb{R}^{m}$, which with respect to the the canonical basis has the representing matrix:}
\[
    Y^*Z =
    \begin{bmatrix}
        \langle y_1,z_1 \rangle & \cdots & \langle y_1,z_\ell \rangle \\
        \vdots                  & \ddots & \vdots                     \\
        \langle y_m,z_1 \rangle & \cdots & \langle y_m,z_\ell \rangle
    \end{bmatrix}.
\]
\rev{Furthermore, given a linear operator $\mathcal{B}: L^2(D) \to L^2(D)$, we write $\mathcal{B} Y $ to denote the quasimatrix whose $i^{\text{th}}$ column is $\mathcal{B} y_i$. }

\rev{Our discussion so far has focused on the case where quasimatrices have a finite number of columns. In general, the definition of a quasimatrix cannot be extended to cases where the quasimatrix has a countably infinite number of columns. However, if the columns form an orthonormal system, it is possible to extend the definition. Consider an orthonormal system $\{u_i\}_{i \in \mathbb{N}} \in L^2(D)$ and let $\ell^2$ denote the set of square-summable sequences. The quasimatrix $U = \begin{bmatrix} u_1 & u_2 & \cdots \end{bmatrix}$ 
is defined to be the bounded linear operator $\ell^2 \to L^2(D)$ satisfying $U\bm{e}_i = u_i$, where $\bm{e}_i$ is the $i^{\text{th}}$ canonical vector in $\ell^2$. Furthermore, for a Hilbert-Schmidt operator $\mathcal{B}:L^2(D) \to L^2(D)$ we define $\mathcal{B}U: \ell^2 \to L^2(D)$ as the linear operator satisfying $\mathcal{B}U \bm{e}_i = \mathcal{B} u_i$. One can verify that $\mathcal{B} U$ is also a Hilbert-Schmidt operator \cite[Def.~4.4.4]{functionaldatanalysis}.}

\subsection{Non-negative self-adjoint trace-class operators}\label{section:spsdoperator}

We consider a non-negative self-adjoint trace-class operator $\mathcal{A} : L^2(D) \to L^2(D)$, i.e., it holds that $\langle \mathcal{A}f,f\rangle\geq 0$ and $\langle \mathcal{A} f,g\rangle = \langle f,\mathcal{A} g\rangle$ for every $f,g\in L^2(D)$,
and the trace norm~\cite[Def.~4.5.1]{functionaldatanalysis} is finite:
\[
    \|\mathcal{A}\|_{\Tr} := \sum_{j=1}^\infty \langle \mathcal{A}e_j,e_j\rangle<\infty,
\]
for any orthonormal basis $\{e_j\}_j$ of $L^2(D)$. Non-negative self-adjoint trace-class operators are Hilbert-Schmidt operators~\cite[Thm.~4.5.2]{functionaldatanalysis}, and therefore admit an eigenvalue decomposition of the form~\cite[Thm.~4.3.1]{functionaldatanalysis}:
\begin{equation} \label{eq:spectraldecomp}
    \mathcal{A} = \sum_{\substack{j=1\\\sigma_j>0}}^\infty \sigma_j\langle u_j,\cdot\rangle u_j,
\end{equation}
where $\sigma_1\geq \sigma_2\geq \cdots\geq 0$ are the eigenvalues of $\mathcal{A}$, and $\{u_j\}_j$ are orthonormal eigenfunctions \rev{forming a complete orthonormal system for the closure of the image of $\mathcal{A}$ \cite[Theorem 4.2.4]{functionaldatanalysis}}. \rev{The series \eqref{eq:spectraldecomp} is converging in the operator norm; see \cite[p.105]{functionaldatanalysis}.} The eigenvalues allow us to express the trace, Hilbert-Schmidt, and operator norms of $\mathcal{A}$ as
\[
    \|\mathcal{A}\|_{\Tr} = \sum_{j=1}^\infty \sigma_j, \quad \|\mathcal{A}\|_{\HS} = \Big(\sum_{j=1}^\infty \sigma_j^2\Big)^{1/2}, \quad \|\mathcal{A}\|_{\op} = \sigma_1,
\]
which are infinite-dimensional analogs of the nuclear, Frobenius, and spectral norms discussed in \Cref{section:finitedimensions}. Furthermore, we introduce the operator $U:\ell^2 \to L^2(D)$ defined by $Uf = \sum_{i=1}^\infty f_i u_i$ for any $f$ in $\ell^2$, the space of square-summable sequences (indexed by positive integers). Then, for a given rank $k\geq 1$, the quasimatrix $U_1:\R^k\to L^2(D)$ contains the first $k$ eigenfunctions of $\mathcal{A}$, and the quasimatrix $U_2:\ell^2\to L^2(D)$ contains the remaining eigenfunctions. \rev{In addition, the adjoint of $U$ is the operator $U^*: L^2(D) \to \ell^2$ defined  $(U^*g)_i = \langle u_i, g\rangle$ for any $g \in L^2(D)$, where $(U^*g)_i$ is the $i^{\text{th}}$ entry in the sequence $U^*g$. $U_1^*$ and $U_2^*$ can be defined in an analogous way.} Finally, we introduce the diagonal matrix and \rev{infinite diagonal matrix} $\bm{\Sigma}_1 = \diag(\sigma_1,\ldots,\sigma_k)$ and $\bm{\Sigma}_2 = \diag(\sigma_{k+1},\sigma_{k+2},\ldots)$, respectively, which contain the eigenvalues of $\mathcal{A}$ in descending order.

\subsection{Gaussian processes} \label{section:gp}

Let $K: D \times D \mapsto \mathbb{R}$ be a continuous symmetric positive semi-definite kernel, with the following eigenvalue decomposition~\cite[Thm.~4.6.5]{functionaldatanalysis}:
\[
    K(x,y) = \sum_{i=1}^{\infty} \lambda_i \psi_i(x) \psi_i(y), \quad \int_D K(x,y) \psi_i(y) dy = \lambda_i \psi_i(x), \quad x,y \in D,
\]
where the sum converges absolutely and uniformly on $D$~\cite{mercer1909functions}. Here, $\lambda_1 \geq \lambda_2 \geq \ldots \geq 0$ are the eigenvalues of the integral operator $\mathcal{K}$ induced by $K$:
\[[\mathcal{K}f](x) = \int_D K(x,y) f(y)\d y, \quad f\in L^2(D),\, x\in D,\]
and $\{\psi_j\}_{j}$ are the corresponding orthonormal eigenfunctions of $\mathcal{K}$ in $L^2(D)$. In the following, we assume that $\mathcal{K}$ is trace-class, i.e., $\sum_{i=1}^{\infty} \lambda_i < \infty$. A Gaussian process (GP) $\omega \sim \mathcal{GP}(0, K)$ is a stochastic process such that for any finite collection of points $x_1,\ldots,x_m \in D$, the vector $\begin{bmatrix} \omega(x_1) & \cdots & \omega(x_m) \end{bmatrix}^\top$ follows a multivariate normal distribution with mean $\bm{0}\in \R^m$ and covariance matrix $\bm{K} = (K(x_i,x_j))_{1\leq i,j\leq m}$. The Karhunen-Lo\`eve expansion~\cite[Ch.3.2]{adler} of $\omega$ is given by
\begin{equation} \label{eq:klexpansion}
    \omega \sim \sum_{k=1}^\infty \lambda_k^{1/2} \zeta_k \psi_k,
\end{equation}
where $\zeta_1,\zeta_2, \ldots \sim \mathcal N(0,1)$ are mutually independent. With probability one, a realization of $\omega$ is in $L^2(D)$ \cite[Thm.~7.2.5]{functionaldatanalysis}. Recalling that $\{u_j\}_j$ denote the eigenfunctions of $\mathcal A$, see~\eqref{eq:spectraldecomp}, we define the \rev{linear operator $\widetilde{\bm{K}}:\ell^2 \to \ell^2$, which under the canonical basis has the following infinite matrix representation}
\[
    \rev{\tK_{i,j}} = \langle u_i, \mathcal{K} u_j\rangle = \sum_{k=1}^\infty \lambda_k \langle \psi_k,u_i\rangle \langle \psi_k,u_j\rangle,\quad i,j\in \N^*.
\]
\rev{Note that $\widetilde{\bm{K}}$ is also a trace-class operator.} We also define the following restrictions of $\tK$ to different sets of indices:
\[
    \tK_{11} = \tK_{\llbracket 1,k\rrbracket\times \llbracket 1,k\rrbracket},\quad \tK_{21} = \tK_{\llbracket k+1,\infty)\times \llbracket 1,k\rrbracket},\quad \tK_{22} = \tK_{\llbracket k+1,\infty)\times \llbracket k+1,\infty)}.
\]
Then, $\tK_{11}$ defines a $k\times k$ matrix, which we assume to be of rank $k$ in the rest of this section. In terms of the quasimatrices $U_1, U_2$ containing the eigenfunctions defined above, we can write
\[
    \tK_{11} = U_1^* \mathcal{K} U_1, \quad \tK_{21} = U_2^* \mathcal{K} U_1, \quad \tK_{22} = U_2^* \mathcal{K} U_2,
\]
in analogy to~\eqref{eq:deftildek}. \rev{Similar to the discussion above, one can define $\mathcal{GP}(0,\widetilde{\bm{K}})$ as follows: a stochastic process $\bm{\omega}$ in $\ell^2$ is said to be distributed as $\mathcal{GP}(0,\widetilde{\bm{K}})$ if, for any finite collection of indices $\ell_1,\ldots,\ell_m \in \mathbb{N}^*$, the vector $\begin{bmatrix} \bm{\omega}_{\ell_1} & \cdots & \bm{\omega}_{\ell_m} \end{bmatrix}$ follows a multivariate normal distribution with mean $\bm{0} \in \mathbb{R}^m$ and covariance matrix $\bm{K} = (\widetilde{\bm{K}}_{\ell_i, \ell_j})_{1 \leq i,j \leq m}$.}\footnote{\rev{By the Kolmogorov two-series theorem, $\bm{\omega}\in \ell^2$ almost surely. }}

\subsection{Infinite-dimensional extension of the Nystr{\"o}m approximation}\label{section:infnystrom}

We are now ready to present the infinite-dimensional extension of the Nystr{\"o}m approximation. Let $k$ be a target rank, $p$ be an oversampling parameter, and $\Omega=\begin{bmatrix}\omega_1&\cdots &\omega_{k+p}\end{bmatrix}$ be a random quasimatrix with $k+p$ columns, whose columns are i.i.d.~from $\mathcal{GP}(0, K)$. The Nystr{\"o}m approximation $\widehat{\mathcal{A}}$ to $\mathcal{A}$ is defined as
\begin{equation} \label{eq_def_nystrom_inf}
    \widehat{\mathcal{A}} := \mathcal{A}\Omega(\Omega^* \mathcal{A}\Omega)^{\dagger}  (\mathcal{A}\Omega)^*.
\end{equation}
Assuming that the realization of $\omega_i$ is in $L^2(D)$, which holds with probability $1$, the Nystr{\"o}m approximation is an operator $\widehat{\mathcal{A}}: L^2(D) \to L^2(D)$ of rank at most $k+p$ \rev{whose action on a function $f \in L^2(D)$ has the explicit representation}
\[
    \rev{[\widehat{\mathcal{A}}f](x) = \sum_{i,j = 1}^{k+p} [\mathcal{A}\omega_i](x) \big[ (\Omega^* \mathcal{A}\Omega)^{\dagger} \big]_{ij} \langle  \omega_j, \mathcal A f \rangle,}
\]
\rev{where $\big[ (\Omega^* \mathcal{A}\Omega)^{\dagger} \big]_{ij}$ is the $(i,j)$-entry of the matrix $(\Omega^* \mathcal{A}\Omega)^{\dagger}$.} As in the finite-dimensional case, the error bounds for the infinite-dimensional analog of the Nystr{\"o}m approximation depend on the prior information of eigenvectors contained in $\mathcal{K}$, which is measured by the following two quantities:
\begin{equation} \label{eq_def_cov_factor_inf}
    \beta_k^{(\xi)} = \frac{\|\bm{\Sigma}_2^{1/2} \tK_{22.1}\bm{\Sigma}_2^{1/2}\|_{\xi}}{\|\bm{\Sigma}_2\|_{\xi}}\|\tK_{11}^{-1}\|_2, \quad \rev{\delta_k^{(\xi)} = \frac{\|\bm{\Sigma}_2^{1/2} \tK_{21}\tK_{11}^{-2}\tK_{21}^*\bm{\Sigma}_2^{1/2}\|_{\xi}}{\|\bm{\Sigma}_2\|_{\xi}}},
\end{equation}
where $\xi \in \{\HS,\Tr,\op\}$.\footnote{\rev{$\beta_k^{(\xi)}$ and $\delta_k^{(\xi)}$ are defined by the Hilbert-Schmidt, trace, and operator norm of trace-class operators $\ell^2 \to \ell^2$. Therefore, these quantities are well-defined.}} These quantities are the infinite-dimensional analogs of \cref{eq_def_cov_factor}.

For each $1\leq j\leq k+p$, we consider the stochastic process $\bm{\omega}_j=\{\langle u_i,\omega_j\rangle,\, i\in \N^*\}$ whose trajectories are in $\ell^2$, and denote by $\bOmega=U^*\Omega=\begin{bmatrix}\bm{\omega}_1&\cdots &\bm{\omega}_{k+p}\end{bmatrix}$ the random \rev{infinite matrix} whose columns are i.i.d. from $\mathcal{GP}(0,\tK)$. Then, we introduce the random $k\times (k+p)$ matrix $\bOmega_1$ as
\[
    \bOmega_1 \coloneqq U_1^* \Omega =
    \begin{bmatrix}
        \langle u_1,\omega_1\rangle & \cdots & \langle u_1,\omega_{k+p}\rangle \\
        \vdots                      & \ddots & \vdots                          \\
        \langle u_k,\omega_1\rangle & \cdots & \langle u_k,\omega_{k+p}\rangle
    \end{bmatrix},
    \quad \omega_j\sim \mathcal{GP}(0,K),\quad 1\leq j\leq k+p,
\]
whose columns are i.i.d.~from $\mathcal{N}(0,\tK_{11})$~\cite[Lem.~1]{boulle2022learning}. Since we assume that $\rank(\tK_{11}) = k$ we know that $\bOmega_1$ has full row-rank with probability one. Therefore, $\bOmega_1$ has almost surely a right inverse $\bOmega_1^{\dagger} = \bOmega_1^*(\bOmega_1\bOmega_1^*)^{-1}$. We define $\bm{\Omega}_2 = U_2^* \Omega$ similarly. Note that the $i^{\text{th}}$ entry of the columns of $\bm{\Omega}_2$ are distributed as $\mathcal{N}(0,(\tK_{22})_{ii})$ and the covariance between the $i^{\text{th}}$ and $j^{\text{th}}$ entry is $(\tK_{22})_{ij}$ \cite[Sec.~3.3]{boulle2022learning}.

\subsubsection{Structural bound}

In this section, we prove an infinite-dimensional analog of the structural bound~\eqref{eq:structurebound}. We begin by stating some basic but useful results on norms of finite sections of $\ell^2$ operators.

\begin{lemma} \label{lemma:operator_limit}
    Let $\bm{B},\bm{C}: \ell^2 \to \ell^2$ be non-negative self-adjoint trace-class operators such that $\bm{B}-\bm{C}$ is non-negative. For $n\in \N^*$, consider the restriction $\bm{B}_{\llbracket 1,n\rrbracket\times \llbracket 1,n\rrbracket}$ of $\bm{B}$ to its first $n$ rows and columns \rev{under the canonical basis of $\ell^2$}. Then, for $\xi \in\{\Tr,\HS,\op\}$, we have
    \begin{subequations}
        \begin{align}
            \lim_{n \rightarrow \infty} \|\bm{B}_{\llbracket 1,n\rrbracket\times \llbracket 1,n\rrbracket}\|_{\xi} & = \|\bm{B}\|_{\xi},\label{eq:3properties_i}                                                                                                                                        \\
            \|\bm{C}\|_{\xi}\leq \|\bm{B}\|_{\xi}                                                                  & \leq \|\bm{B}_{\llbracket 1,n\rrbracket\times \llbracket 1,n\rrbracket}\|_{\xi} + \|\bm{B}_{\llbracket n+1,\infty)\times \llbracket n+1,\infty)}\|_{\xi}.\label{eq:3properties_ii}
        \end{align}
    \end{subequations}
\end{lemma}
\begin{proof}
For $\xi \in \{\HS,\Tr\}$, the property~\eqref{eq:3properties_i}  follows from the absolute convergence of the involved series \rev{$\lim_{n \rightarrow\infty} \sum_{i,j = 1}^{n} \bm{B}_{ij}^2 = \sum_{i,j = 1}^{\infty} \bm{B}_{ij}^2$ and $\lim_{n \rightarrow\infty} \sum_{i = 1}^{n} \bm{B}_{ii} = \sum_{i = 1}^{\infty} \bm{B}_{ii}$, where the entries $\bm{B}_{ij}$ are defined with respect to the canonical basis of $\ell^2$.} For $\xi = \op$, \rev{let $\bm{P}_n$ be the orthogonal projection onto the first $n$ canonical basis vectors in $\ell^2$. Then, $\|\bm{B}_{\llbracket 1,n\rrbracket\times \llbracket 1,n\rrbracket}\|_{\op} = \|\bm{P}_n \bm{B} \bm{P}_n\|_{\op}$. Hence,} the triangle inequality yields\rev{
\begin{align*}
    |\|\bm{B}\|_{\op}-\|\bm{B}_{\llbracket 1,n\rrbracket\times \llbracket 1,n\rrbracket}\|_{\op}| &\leq \|\bm{B} - \bm{P}_n\bm{B}\bm{P}_n\|_{\op} \\
    & \leq \|\bm{B} - \bm{P}_n\bm{B}\bm{P}_n\|_{\HS}\\
    & = \sqrt{\|\bm{B}\|_{\HS}^2 - \|\bm{B}_{\llbracket 1,n\rrbracket\times \llbracket 1,n\rrbracket}\|_{\HS}^2}\rightarrow 0,\quad \text{as }n\to\infty.
\end{align*}}
\rev{To obtain \eqref{eq:3properties_ii}, first pick some $m \geq n$. Combining \cite[Thm.~2.1]{lee} and \cite[Lem.~3.1]{persson2022randomized}, we have $$\|\bm{C}_{\llbracket 1,m\rrbracket\times \llbracket 1,m\rrbracket}\|_{\xi}\leq \|\bm{B}_{\llbracket 1,n\rrbracket\times \llbracket 1,n\rrbracket}\|_{\xi} \leq \|\bm{B}_{\llbracket 1,n\rrbracket\times \llbracket 1,n\rrbracket}\|_{\xi} + \|\bm{B}_{\llbracket n+1,m)\times \llbracket n+1,m)}\|_{\xi}.$$
Taking the limit as $m \to \infty$ and using \eqref{eq:3properties_i} yields the desired result. }
\end{proof}
We can now proceed with an infinite-dimensional analog of the finite-dimensional structural bound \eqref{eq:structurebound}.

\begin{lemma}[Infinite-dimensional structural bound]\label{lemma:infnystrom_structural}
    Let $\mathcal{A}:L^2(D)\to L^2(D)$ be a self-adjoint non-negative trace-class operator, $k,p\geq 1$, and $\Omega:\R^{k+p}\to L^2(D)$ be a quasimatrix with $k+p$ columns such that the matrix $\bOmega_1=U_1^*\Omega\in \R^{k\times (k+p)}$ is full rank. Then,
    \[
        \|\mathcal{A}- \mathcal{A} \Omega(\Omega^*\mathcal{A} \Omega)^{\dagger}(\mathcal{A}\Omega)^*\|_{\xi} \leq \|\bm{\Sigma}_2\|_{\xi} + \|(\bm{\Sigma}_2^{1/2} \bOmega_2 \bOmega_1^{\dagger})^*\bm{\Sigma}_2^{1/2} \bOmega_2 \bOmega_1^{\dagger}\|_{\xi},
    \]
    where $\xi\in \{\mathrm{op}, \mathrm{HS},\mathrm{Tr}\}$.
\end{lemma}

\begin{proof}
    Let $\bOmega = U^* \Omega$ be defined as in \Cref{section:infnystrom} and $\mathcal{P}_{U}=UU^*:L^2(D)\to L^2(D)$ denote the orthogonal projection onto the range of $\mathcal{A}$. Since $\mathcal{A}$ is a self-adjoint operator, we have
    $
        \mathcal{A} = \mathcal{P}_{U} \mathcal{A} = \mathcal{A} \mathcal{P}_{U} = \mathcal{P}_{U} \mathcal{A} \mathcal{P}_{U}
    $. Therefore,
    \begin{align*}
         & \|\mathcal{A}- \mathcal{A} \Omega(\Omega^*\mathcal{A} \Omega)^{\dagger}(\mathcal{A}\Omega)^*\|_{\xi} = \|\mathcal{P}_{U}\mathcal{A}\mathcal{P}_{U}- \mathcal{P}_{U}\mathcal{A}\mathcal{P}_{U} \Omega(\Omega^*\mathcal{P}_{U}\mathcal{A}\mathcal{P}_{U} \Omega)^{\dagger}(\mathcal{P}_{U}\mathcal{A}\mathcal{P}_{U}\Omega)^*\|_{\xi} \\
         & \quad\qquad = \|U(U^* \mathcal{A} U - U^* \mathcal{A} U \bOmega(\bOmega^* U^* \mathcal{A} U \bOmega)^{\dagger}(U^* \mathcal{A} U \bOmega)^* )U^*\|_{\xi}                                                                                                                                                                            \\
         & \quad\qquad = \|U^* \mathcal{A} U - U^* \mathcal{A} U \bOmega(\bOmega^* U^* \mathcal{A} U \bOmega)^{\dagger}(U^*  \mathcal{A} U \bOmega)^* \|_{\xi} = \|\bm{\Sigma} - \bm{\Sigma} \bOmega(\bOmega^* \bm{\Sigma}\bOmega)^{\dagger}(\bm{\Sigma}\bOmega)^*\|_{\xi},
    \end{align*}
    where the third equality follows from the unitary invariance of the norms
    and the fourth equality is due to the relation $U^* \mathcal{A} U = \bm{\Sigma}$. \rev{Now note that $\|\bm{\Sigma} - \bm{\Sigma} \bOmega(\bOmega^* \bm{\Sigma}\bOmega)^{\dagger}(\bm{\Sigma}\bOmega)^*\|_{\xi} = \|\bm{\Sigma}^{1/2}(\bm{I} - \bm{P}_{\bm{\Sigma}^{1/2} \bm{\Omega}}) \bm{\Sigma}^{1/2}\|_{\xi},$ where $\bm{I}$ is the identity operator and $\bm{P}_{\bm{\Sigma}^{1/2} \bm{\Omega}}$ is the orthogonal projection onto $\range(\bm{\Sigma}^{1/2} \bm{\Omega})$. Let $$\bm{Z} = \bm{\Sigma} \bm{\Omega} \bm{\Omega}_1^{\dagger} \bm{\Sigma}_1^{-1/2} = \begin{bmatrix} \bm{I} \\ \bm{F} \end{bmatrix}, \quad \bm{F} = \bm{\Sigma}_2^{1/2} \bm{\Omega}_2  \bm{\Omega}_1^{\dagger} \bm{\Sigma}_1^{-1/2},$$ 
    where we used that $\bm{\Omega}_1$ is full rank. We have $\range(\bm{Z}) \subseteq \range(\bm{\Sigma}^{1/2} \bm{\Omega})$ and therefore $\bm{I} - \bm{P}_{\bm{Z}} \succeq \bm{I} - \bm{P}_{\bm{\Sigma}^{1/2} \bm{\Omega}}$, where $\bm{P}_{\bm{Z}} = \bm{Z}(\bm{Z}^* \bm{Z})^{\dagger} \bm{Z}^*$ is the orthogonal projection onto $\range(\bm{Z})$. Hence, by \Cref{lemma:operator_limit} $\|\bm{\Sigma}^{1/2}(\bm{I} - \bm{P}_{\bm{\Sigma}^{1/2} \bm{\Omega}}) \bm{\Sigma}^{1/2}\|_{\xi} \leq \|\bm{\Sigma}^{1/2}(\bm{I} - \bm{P}_{\bm{Z}}) \bm{\Sigma}^{1/2}\|_{\xi}$.} As in the finite-dimensional case, the rest of the proof follows the argument of the proof of~\cite[Lem.~3.13]{persson2022randomized} for the operator monotone function $f:x\mapsto x$ and $q = 1$, now using \rev{the $2 \times 2$ block-structure of $\bm{I}-\bm{P}_{\bm{Z}}$: $$\bm{I}-\bm{P}_{\bm{Z}} = \begin{bmatrix} \bm{I} - (\bm{I} + \bm{F}^* \bm{F})^{-1} & -(\bm{I} + \bm{F}^* \bm{F})^{-1} \bm{F}^* \\ -\bm{F}(\bm{I} + \bm{F}^* \bm{F})^{-1} & \bm{I} - \bm{F}(\bm{I} + \bm{F}^* \bm{F})^{-1}\bm{F}^* \end{bmatrix}$$ and} the fact that the inequalities in \eqref{eq:3properties_ii} are infinite-dimensional analogs of~\cite[Lems.~A.1 and 3.1]{persson2022randomized} with $f:x\mapsto x$ and $\|\cdot\| \coloneqq \|\cdot\|_{\xi}$.
\end{proof}

\subsection{Probabilistic bounds}
With the structural bound in place, we proceed to derive probabilistic bounds for the infinite-dimensional Nystr{\"o}m approximation~\eqref{eq_def_nystrom_inf}.

\begin{theorem}[Infinite-dimensional Nystr{\"o}m approximation]\label{theorem:inf_probabilistic}
    Let $\mathcal{A}: L^2(D) \to L^2(D)$ be a non-negative self-adjoint trace-class operator, $2\leq k\leq \rank(\mathcal A)$ be a target rank, and $p\ge 4$ be an oversampling parameter. Let $\Omega$ be a quasimatrix with $k+p$ columns  i.i.d.~from $\mathcal{GP}(0,K)$, with a kernel $K$ such that the matrix $\tK_{11}$ defined in \eqref{eq:deftildek} is invertible. Then, the Nystr{\"o}m approximation $\widehat{\mathcal{A}}:L^2(D)\to L^2(D)$ defined in \cref{eq_def_nystrom_inf} satisfies
    \begin{subequations} \label{eq_expectation_tail bound_inf}
        \begin{align}
            \mathbb{E}[\|\mathcal{A}-\widehat{\mathcal{A}}\|_\op] & \leq \left(1 + \frac{3k}{p-1} \beta_k^{(\op)} + 3\delta_k^{(\op)}\right)\|\bm{\Sigma}_2\|_\op + \frac{3e^2(k+p)}{p^2-1}\beta_k^{(\Tr)}\|\bm{\Sigma}_2\|_\Tr,\label{eq:inf_expectation_op} \\
            \mathbb{E}[\|\mathcal{A}-\widehat{\mathcal{A}}\|_\Tr] & \leq \left(1 + \frac{k}{p-1} \beta_k^{(\Tr)} + \delta_k^{(\Tr)}\right)\|\bm{\Sigma}_2\|_\Tr,\label{eq:inf_expectation_tr}\\
            \rev{\mathbb{E}[\|\mathcal{A}-\widehat{\mathcal{A}}\|_\HS]} & \rev{\leq \left( 1+2 \delta_k^{(\HS)} + 2\sqrt{c_1} \beta_k^{(\HS)}\right) \|\bm\Sigma_2\|_\HS+2\sqrt{c_2} \beta_k^{(\Tr)}\|\bm\Sigma_2\|_\Tr,\label{eq:inf_expectation_hs}}
        \end{align}
    \end{subequations}
    where $c_1=\mathcal{O}(k^2/p^2)$, $c_2=\mathcal{O}(k^2/p^2)$ are the constants defined in \cref{eq_const_c1_c2}. Let $u,t\geq 1$, then
    \begin{subequations} \label{eq_tail bound_inf}
        \begin{align}
             \|\mathcal{A}-\widehat{\mathcal{A}}\|_\op & \leq \left(1+4\delta_k^{(\op)}+4(d_1+d_2 u^2) t^2\beta_k^{(\op)}\right)\|\bm{\Sigma}_2\|_\op+4 d_2 t^2\beta_k^{(\Tr)}\|\bm{\Sigma}_2\|_\Tr,\label{eq:inf_tail bound_op} \\
            \|\mathcal{A}-\widehat{\mathcal{A}}\|_\Tr & \leq \left(1 + 2 \delta_k^{(\Tr)} + d_1 t^2\beta_k^{(\Tr)}\right)\|\bm{\Sigma}_2\|_\Tr+2 d_2 t^2 u^2 \beta_k^{(\op)}\|\bm{\Sigma}_2\|_\op,\label{eq:inf_tail bound_tr}\\
            \rev{\|\mathcal{A}-\widehat{\mathcal{A}}\|_\HS} & \rev{\leq
            \|\bm{\Sigma}_2\|_\HS + 4 \left(\delta_k^{(\HS)} + t^2(d_1+d_3)\beta_k^{(\HS)}\right)\|\bm\Sigma_2\|_\HS + 4 t^2 d_3\beta_k^{(\Tr)}\|\bm\Sigma_2\|_\Tr} \label{eq:inf_tail bound_hs}\\
            & \rev{\quad + 2 t^2 u^2 d_2\beta_k^{(2)}\|\bm\Sigma_2\|_\op,} \nonumber      
        \end{align}
    \end{subequations}
    with probability $\geq 1-3t^{-p}-e^{-u^2/2}$. Here, $d_1 = \mathcal{O}(k/p)$, $d_2 = \mathcal{O}(k/p)$, and $d_3 = \mathcal{O}(k^{3/2}/p)$ are the constants defined in \cref{eq_def_di} \rev{and in \Cref{theorem:tail bound_frob}}.
\end{theorem}

The proof of \cref{theorem:inf_probabilistic} occupies the rest of this section and follows from a continuity argument on the generalization of the Nystr{\"o}m approximation to correlated Gaussian vectors analyzed earlier in \Cref{section:finitedimensions}. Let $n\geq 1$ and $s\in[1,\infty]$, we first define the following random variables:
\begin{equation}\label{eq:rvs}
    X_{s,n} = \|(\bm{\Sigma}_2^{1/2})_{\llbracket 1,n\rrbracket\times \llbracket 1,n\rrbracket} (\bm{\Omega}_2)_{\llbracket 1,n\rrbracket\times \llbracket 1,k+p\rrbracket} \bm{\Omega}_1^{\dagger}\|_{(2s)}, \quad X_{s} = \|\bm{\Sigma}_2^{1/2} \bm{\Omega}_2 \bm{\Omega}_1^{\dagger}\|_{(2s)},
\end{equation}
where $\|\cdot\|_{(2s)}$ denotes the Schatten-$2s$ norm. We aim to show the convergence of $X_{s,n}$ to $X_{s}$ as $n\to\infty$, and begin with a preliminary result on the finiteness of the expectation of $X_{s}$.

\begin{lemma}[Expectation of $X_s$]\label{lemma:expectationbounded}
    For $s\in[1,\infty]$, let $X_s$ be the random variable defined in \eqref{eq:rvs}. Then, if $p \geq 4$, we have $\mathbb{E}_\Omega^2[X_{s}] < \infty$.
\end{lemma}
\begin{proof}
    We first notice that $2s \geq 2$ implies
    \[
        \|\bm{\Sigma}_2^{1/2}\bOmega_2\bOmega_1^{\dagger}\|_{(2s)}^2 \leq \|\bm{\Sigma}_2^{1/2}\|_{\op}^2 \|\bOmega_2\|_{\HS}^2\|\bOmega_1^{\dagger}\|_{\Frob}^2.
    \]
    Noting that $\bOmega_1$, $\bOmega_2$ are not independent, we need to establish $\mathbb{E}[\|\bOmega_2\|_{\HS}^4] < \infty$ and $\mathbb{E}[\|\bOmega_1^{\dagger}\|_{\Frob}^4] < \infty$ in order to conclude the result from H{\"o}lder's inequality. First, \cref{lemma:inversewishart_norm_expectation} ensures that $\mathbb{E}[\|\bOmega_1^{\dagger}\|_{\Frob}^4] < \infty$ since $\bOmega_1$ is a $k\times (k+p)$ matrix whose columns are i.i.d.~from $\mathcal{N}(\bm{0}, \tK_{11})$ \cite[Lem.~1]{boulle2022learning}. By the Karhunen–Lo\`eve expansion~\eqref{eq:klexpansion}, an arbitrary column $\omega$
    of $\Omega$ satisfies
    \begin{align*}
        \mathbb{E}\big[\|\omega\|^4_{L^2(D)}\big] & = \sum_{i,j = 1}^\infty \mathbb{E}[\zeta^2_i \zeta^2_j] \lambda_i \lambda_j
        = \sum_{i=1}^\infty \mathbb{E}[\zeta^4_i] \lambda_i^2 +  \sum_{i \not= j} \mathbb{E}[\zeta_i^2 \zeta_j^2]\lambda_i \lambda_j                                   \\
                                                  & = 3 \sum_{i=1}^\infty \lambda_i^2 + \sum_{i \not= j} \lambda_i \lambda_j \le 3 \| \mathcal{K} \|_{\Tr}^2 < \infty.
    \end{align*}
    In turn,
    \[
        \mathbb{E}\|\bm{\Omega}_2\|_{\HS}^4 \leq \mathbb{E}\|\Omega\|_{\HS}^4 \leq \mathbb{E}\big( \|\omega_1\|^2_{L^2(D)} + \cdots + \|\omega_{k+p} \|^2_{L^2(D)} \big)^2 \le  (k+p)^2 \mathbb{E}\big[\|\omega\|^4_{L^2(D)}\big] < \infty.
    \]
\end{proof}

We proceed by showing that $\lim_{n \rightarrow \infty}\mathbb{E}^2[X_{s} - X_{s,n}] = 0$.

\begin{lemma}[Convergence of $X_{s,n}$ to $X_s$]\label{lemma:limits}
    For $s\in[1,\infty]$, let $X_s$, $X_{s,n}$ be the random variables defined in \eqref{eq:rvs} for $n\geq 1$. Then, if $p \geq 4$, we have $\lim_{n \rightarrow \infty}\mathbb{E}^2[X_{s} - X_{s,n}] = 0$.
\end{lemma}

\begin{proof}
    \rev{For $n\geq 1$, we define the matrix $\bOmega_2^{(n)}$, with an infinite number of rows, whose first $n$ rows are equal to the first $n$ rows of $\bOmega_2$ and the remaining rows are zero.} Then,
    \[
        X_{s,n} = \|(\bm{\Sigma}_2^{1/2})_{\llbracket 1,n\rrbracket\times\llbracket 1,n\rrbracket}(\bOmega_2)_{\llbracket 1,n\rrbracket\times\llbracket 1,k+p\rrbracket} \bOmega_1^{\dagger}\|_{(2s)} = \|\bm{\Sigma}_2^{1/2}\bOmega_2^{(n)} \bOmega_1^{\dagger}\|_{(2s)}.
    \]
    Combining the triangle inequality and sub-multiplicativity of the Schatten-$s$ norm, and using the fact that $2s \geq 2$, we have
    \begin{equation*}
        (X_{s}-X_{s,n})^2 \leq \|\bm{\Sigma}_2^{1/2}\|_{\op}^2 \|\bOmega_1^{\dagger}\|_{\Frob}^2 \|\bOmega_2-\bOmega_2^{(n)}\|_{\HS}^2.
    \end{equation*}
    By Hölder's inequality it suffices to show that $\lim_{n \rightarrow \infty}\mathbb{E}[\|\bOmega_2 - \bOmega_2^{(n)}\|_{\HS}^4] = 0$, since $\mathbb{E}[\|\bOmega_1^{\dagger}\|_\Frob^4] < \infty$ by \cref{lemma:inversewishart_norm_expectation}. Let $\bm{\omega}_i$ and $\bm{\omega}_i^{(n)}$ denote the $i^{\text{th}}$ columns of $\bOmega_2$ and $\bOmega_2^{(n)}$, respectively. Using the monotonicity and triangle inequality of $L^p$-norms, we have
    \begin{align*}
        \mathbb{E}^4[\|\bm{\Omega}_2 - \bm{\Omega}_2^{(n)}\|_{\HS}] & \leq \mathbb{E}^4\left[\sum_{i=1}^{k+p} \|\bm{\omega}_i - \bm{\omega}_i^{(n)}\|_2\right]\leq \sum_{i=1}^{k+p} \mathbb{E}^4\left[\|\bm{\omega}_i - \bm{\omega}_i^{(n)}\|_2\right] \\
                                                                    & = (k+p) \mathbb{E}^4\left[\|\bm{\omega}_1 - \bm{\omega}_1^{(n)}\|_2\right],
    \end{align*}
    since the columns of $\bm{\Omega}_2$ are identically distributed.

    We are now going to verify that
    \[
        \lim_{n \rightarrow \infty}\mathbb{E}\left[\|\bm{\omega}_1 - \bm{\omega}_1^{(n)}\|_2^4\right] = 0.
    \]
    Following~\cite[Sec.~3.3]{boulle2022learning}, we know that the entries of $\bm{\omega}_1$ satisfy $(\bm{\omega}_1)_i \sim \mathcal{N}(0,(\tK_{22})_{ii})$ for $i\in\N^*$. Let $Y = \sum_{i=n+1}^{\infty} (\bm{\omega}_1)_i^2 = \|\bm{\omega}_1 - \bm{\omega}_1^{(n)}\|_2^2$. Combining the non-negativity of the summands and the Fubini-Tonelli theorem, we can interchange the summation and expectation to obtain
    \begin{align*}
        \mathbb{E}\left[\|\bm{\omega}_1 - \bm{\omega}_1^{(n)}\|_2^4\right] & = \mathbb{E} [Y^2] = \mathbb{E}\left[\sum_{i=n+1}^{\infty} (\bm{\omega}_1)_i^2 Y\right] = \sum_{i=n+1}^{\infty} \mathbb{E}\left[(\bm{\omega}_1)_i^2 Y\right] \\
                                                                           & \leq \sum_{i=n+1}^{\infty} \sqrt{\mathbb{E}[(\bm{\omega}_1)_i^4]} \sqrt{\mathbb{E} [Y^2]},
    \end{align*}
    where the last inequality follows from the Cauchy-Schwarz inequality. Hence,
    \begin{equation*}
        \mathbb{E}\left[\|\bm{\omega}_1 - \bm{\omega}_1^{(n)}\|_2^4\right] \leq \left(\sum_{i=n+1}^{\infty} \sqrt{\mathbb{E}[(\bm{\omega}_1)_i^4]}\right)^2 = 3 \left(\sum_{i=n+1}^{\infty} (\tK_{22})_{ii}\right)^2 \rightarrow 0,\quad \text{as } n \rightarrow \infty,
    \end{equation*}
    since $\tr(\tK_{22})\leq \tr(\mathcal{K})<\infty$.
\end{proof}

Combining~\cref{lemma:expectationbounded,lemma:limits}, we obtain that $\lim_{n \rightarrow \infty}\mathbb{E}^2[X_{s,n}] = \mathbb{E}^2[X_{s}]$. Hence, if we have a family of bounds $\mathbb{E}^2[X_{s,n}] \leq y_{s,n}$ with $y_{s,n}\to y_{s}$ as $n\to\infty$, then $\mathbb{E}^2[X_{s}] \leq y_{s}$. Furthermore, combining the continuous mapping theorem and the fact that \rev{$L_2$-}convergence implies convergence in distribution, for any positive sequence $z_{s,n}\to z_s < \infty$, we have $\lim_{n \rightarrow \infty}\mathbb{P}(X_{s,n}^2 > z_{s,n}) = \mathbb{P}(X_{s}^2 > z_{s})$. We can now proceed with the proof of \cref{theorem:inf_probabilistic}, which uses results from \Cref{section:finitedimensions} to derive expressions for $y_{s,n}$ and $z_{s,n}$ and show that, for $s \in \{2,\infty,1\}$, they converge to the right-hand sides of \cref{eq_expectation_tail bound_inf,eq_tail bound_inf}.

\begin{proof}[Proof of \Cref{theorem:inf_probabilistic}]
    Let $\|\bm{\Sigma}_2\|_{\HS} + y_{2}$, $\|\bm{\Sigma}_2\|_{\op} + y_{\infty}$, and $\|\bm{\Sigma}_2\|_{\Tr} + y_{1}$ be the respective right-hand sides in \cref{eq_expectation_tail bound_inf}, and $\|\bm{\Sigma}_2\|_{\HS} + z_{2}$, $\|\bm{\Sigma}_2\|_{\op} + z_{\infty}$, and $\|\bm{\Sigma}_2\|_{\Tr} + z_{1}$ be the right-hand sides in \cref{eq_tail bound_inf}. For $s \in \{2,\infty,1\}$, we aim to prove that $\mathbb{E}[X_{s}^2] \leq y_s$ and $\mathbb{P}(X_s^2 \geq z_s) \leq 3t^{-p} + e^{-u^2/2}$. Following \cref{theorem:expectation_frob}, we have
    \begin{align*}
        \mathbb{E}[X_{2,n}^2] & \leq 2 \|(\bm{\Sigma}_2^{1/2})_{\llbracket 1,n\rrbracket\times \llbracket 1,n\rrbracket} (\widetilde{\bm{K}}_{21})_{\llbracket 1,n\rrbracket\times \llbracket 1,k\rrbracket} \rev{\widetilde{\bm{K}}_{11}^{-2}}(\widetilde{\bm{K}}_{21})_{\llbracket 1,n\rrbracket\times \llbracket 1,k\rrbracket}^* (\bm{\Sigma}_2^{1/2})_{\llbracket 1,n\rrbracket\times \llbracket 1,n\rrbracket}\|_{\Frob} \\
                              & \quad + 2\sqrt{c_1} \|(\bm{\Sigma}_2^{1/2})_{\llbracket 1,n\rrbracket\times \llbracket 1,n\rrbracket} (\widetilde{\bm{K}}_{22.1})_{\llbracket 1,n\rrbracket\times \llbracket 1,n\rrbracket} (\bm{\Sigma}_2^{1/2})_{\llbracket 1,n\rrbracket\times \llbracket 1,n\rrbracket}\|_{\Frob}\rev{\|\widetilde{\bm{K}}_{11}^{-1}\|_2}                                                                                                   \\
                              & \quad + 2\sqrt{c_2} \|(\bm{\Sigma}_2^{1/2})_{\llbracket 1,n\rrbracket\times \llbracket 1,n\rrbracket} (\widetilde{\bm{K}}_{22.1})_{\llbracket 1,n\rrbracket\times \llbracket 1,n\rrbracket}(\bm{\Sigma}_2^{1/2})_{\llbracket 1,n\rrbracket\times \llbracket 1,n\rrbracket}\|_{*}\rev{\|\widetilde{\bm{K}}_{11}^{-1}\|_2} := y_{2,n},
    \end{align*}
    and with probability greater than $1-3^{-p} - e^{-u^2/2}$, we have
    \begin{align*}
        X_{2,n}^2 & \geq 4 \|(\bm{\Sigma}_2^{1/2})_{\llbracket 1,n\rrbracket\times \llbracket 1,n\rrbracket} (\widetilde{\bm{K}}_{21})_{\llbracket 1,n\rrbracket\times \llbracket 1,k\rrbracket} \rev{\widetilde{\bm{K}}_{11}^{-2}}(\widetilde{\bm{K}}_{21})_{\llbracket 1,n\rrbracket\times \llbracket 1,k\rrbracket}^* (\bm{\Sigma}_2^{1/2})_{\llbracket 1,n\rrbracket\times \llbracket 1,n\rrbracket}\|_{\Frob} \\
                  & \quad + 4t^2(d_1+d_3)\|(\bm{\Sigma}_2^{1/2})_{\llbracket 1,n\rrbracket\times \llbracket 1,n\rrbracket} (\widetilde{\bm{K}}_{22.1})_{\llbracket 1,n\rrbracket\times \llbracket 1,n\rrbracket} (\bm{\Sigma}_2^{1/2})_{\llbracket 1,n\rrbracket\times \llbracket 1,n\rrbracket}\|_{\Frob}\rev{\|\widetilde{\bm{K}}_{11}^{-1}\|_2}                                                                                                   \\
                  & \quad + 4 t^2 d_3\|(\bm{\Sigma}_2^{1/2})_{\llbracket 1,n\rrbracket\times \llbracket 1,n\rrbracket} (\widetilde{\bm{K}}_{22.1})_{\llbracket 1,n\rrbracket\times \llbracket 1,n\rrbracket} (\bm{\Sigma}_2^{1/2})_{\llbracket 1,n\rrbracket\times \llbracket 1,n\rrbracket}\|_{*}\rev{\|\widetilde{\bm{K}}_{11}^{-1}\|_2}                                                                                                           \\
                  & \quad + 2 t^2 u^2 d_2\|(\bm{\Sigma}_2^{1/2})_{\llbracket 1,n\rrbracket\times \llbracket 1,n\rrbracket} (\widetilde{\bm{K}}_{22.1})_{\llbracket 1,n\rrbracket\times \llbracket 1,n\rrbracket} (\bm{\Sigma}_2^{1/2})_{\llbracket 1,n\rrbracket\times \llbracket 1,n\rrbracket}\|_{2}\rev{\|\widetilde{\bm{K}}_{11}^{-1}\|_2} := z_{2,n}.
    \end{align*}
    Following \cref{theorem:expectation_nystrom,theorem:tail_nystrom} we know that
    \begin{align*}
         & \mathbb{E}[X_{\infty,n}] \leq y_{\infty,n},                                  &  & \mathbb{E}[X_{1,n}] \leq y_{1,n},                                  \\
         & \mathbb{P}\left(X_{\infty,n} > z_{\infty,n}\right) \leq 3^{-p} + e^{-u^2/2}, &  & \mathbb{P}\left(X_{1,n} > z_{1,n}\right) \leq 3^{-p} + e^{-u^2/2},
    \end{align*}
    where $y_{\infty,n}$, $y_{1,n}$, $z_{\infty,n}$, and $z_{1,n}$ can be defined analogously to $y_{2,n}$ and $z_{2,n}$ using \cref{eq_exp_spect_nucl,eq_tail_spect_nucl}. Moreover, \cref{lemma:limits} implies that $\lim_{n \rightarrow \infty} \mathbb{E}[X_{s,n}^2] = \mathbb{E}[X_s^2]$ and convergence in distribution, which implies $\lim_{n \rightarrow \infty}\mathbb{P}(X_{s,n} > z_{s,n}) = \mathbb{P}(X_{s} > \lim_{n \rightarrow \infty}z_{s,n})$. Hence, it is sufficient to show that $\lim_{n \rightarrow \infty}y_{s,n} = y_{s}$ and $\lim_{n \rightarrow \infty}z_{s,n} = z_s$ for $s \in \{2,\infty,1\}$. For this purpose, it is sufficient to show that
    \begin{align}
         & \lim_{n \rightarrow \infty}\|(\bm{\Sigma}_2^{1/2})_{\llbracket 1,n\rrbracket\times \llbracket 1,n\rrbracket} (\widetilde{\bm{K}}_{22.1})_{\llbracket 1,n\rrbracket\times \llbracket 1,n\rrbracket} (\bm{\Sigma}_2^{1/2})_{\llbracket 1,n\rrbracket\times \llbracket 1,n\rrbracket}\|_{\xi}                                                              = \|\bm{\Sigma}_2^{1/2}\widetilde{\bm{K}}_{22.1}\bm{\Sigma}_2^{1/2}\|_{\xi}, \nonumber \\
         & \lim_{n \rightarrow \infty}\|(\bm{\Sigma}_2^{1/2})_{\llbracket 1,n\rrbracket\times \llbracket 1,n\rrbracket} (\widetilde{\bm{K}}_{21})_{\llbracket 1,n\rrbracket\times \llbracket 1,k\rrbracket} \rev{\widetilde{\bm{K}}_{11}^{-2}}(\widetilde{\bm{K}}_{21})_{\llbracket 1,n\rrbracket\times \llbracket 1,k\rrbracket}^* (\bm{\Sigma}_2^{1/2})_{\llbracket 1,n\rrbracket\times \llbracket 1,n\rrbracket}\|_{\xi}  = \label{eq:split}                 \\
         & \hspace{0.63\textwidth}\|\bm{\Sigma}_2^{1/2}\widetilde{\bm{K}}_{21}\rev{\widetilde{\bm{K}}_{11}^{-2}}\widetilde{\bm{K}}_{21}^* \bm{\Sigma}_2^{1/2}\|_{\xi}, \nonumber
    \end{align}
    for $\xi \in \{\HS,\Tr,\op\}$, as inserting the definitions of $\beta_k^{(\xi)}$ and $\delta_k^{(\xi)}$ would then give the result. Note that since $\bm{\Sigma}_2$ is diagonal we have
    \begin{align*}
         & (\bm{\Sigma}_2^{1/2})_{\llbracket 1,n\rrbracket\times \llbracket 1,n\rrbracket} (\widetilde{\bm{K}}_{22.1})_{\llbracket 1,n\rrbracket\times \llbracket 1,n\rrbracket} (\bm{\Sigma}_2^{1/2})_{\llbracket 1,n\rrbracket\times \llbracket 1,n\rrbracket} = (\bm{\Sigma}_2^{1/2}\widetilde{\bm{K}}_{22.1}\bm{\Sigma}_2^{1/2})_{\llbracket 1,n\rrbracket\times \llbracket 1,n\rrbracket}, \\
         & (\bm{\Sigma}_2^{1/2})_{\llbracket 1,n\rrbracket\times \llbracket 1,n\rrbracket} (\widetilde{\bm{K}}_{21})_{\llbracket 1,n\rrbracket\times \llbracket 1,k\rrbracket} \widetilde{\bm{K}}_{11}^{-1}(\widetilde{\bm{K}}_{21})_{\llbracket 1,n\rrbracket\times \llbracket 1,k\rrbracket}^* (\bm{\Sigma}_2^{1/2})_{\llbracket 1,n\rrbracket\times \llbracket 1,n\rrbracket} =              \\
         & \hspace{0.6\textwidth} (\bm{\Sigma}_2^{1/2}\widetilde{\bm{K}}_{21} \rev{\widetilde{\bm{K}}_{11}^{-2}}\widetilde{\bm{K}}_{21}^* \bm{\Sigma}_2^{1/2})_{\llbracket 1,n\rrbracket\times \llbracket 1,n\rrbracket}.
    \end{align*}
    Since $\bm{\Sigma}_2^{1/2}\widetilde{\bm{K}}_{22.1}\bm{\Sigma}_2^{1/2}$ and $\bm{\Sigma}_2^{1/2}\widetilde{\bm{K}}_{21} \rev{\widetilde{\bm{K}}_{11}^{-2}}\widetilde{\bm{K}}_{21}^* \bm{\Sigma}_2^{1/2}$ are non-negative trace-class operators, applying \eqref{eq:3properties_i} yields \eqref{eq:split}, as desired.
\end{proof}

\section{Numerical experiments}\label{section:experiments}

In this section, we test the infinite-dimensional Nystr{\"o}m approximation proposed in this work to compute low-rank approximations to different \rev{trace class} operators. \cref{alg:nystrom} presents the pseudocode of a variant of~\cite[Alg.~16]{randnla} for non-negative self-adjoint trace-class operators. We emphasize that the orthonormalization steps in line~\ref{line:orthonormalization} and the shift in line~\ref{line:shift} improve the numerical stability of the algorithm; see e.g. \cite{tropp2017fixed}. Moreover, our implementation of the Nystr{\"o}m approximation computes the Cholesky factorization of the symmetric part of $\bm{\Omega}^* \bm{Y}_{\nu}$ in line~\ref{line:chol}. Finally, we note that \rev{in the ideal setting when we are able to compute with functions exactly}, with $\nu = 0$ in line~\ref{line:nu}, the approximation returned by \cref{alg:nystrom} is mathematically equivalent to \cref{eq_def_nystrom_inf}.

\rev{The aim of our numerical experiments is to validate \cref{alg:nystrom}. For this purpose, we consider a set $D = [-1,1]^d$ for $d = 1,2$ and compute low-rank approximations to different operators $\mathcal{A}: L^2(D) \to L^2(D)$. In most of the experiments, we will carry out all operations on functions on $D$ (approximately) using the Chebfun software package \cite{driscoll}, which is an open-source MATLAB package containing algorithms for performing numerical linear algebra operations on functions and quasimatrices, such as the QR decomposition~\cite{townsend2013extension,quasimatrix}. A common choice for the covariance kernel associated with the random field $\Omega$ in \cref{alg:nystrom} is the squared-exponential kernel:
\begin{equation}\label{eq:sekernel}
    K(\bm{x},\bm{y}) = \prod\limits_{i=1}^d \exp\left(-\frac{(\bm{x}_i-\bm{y}_i)^2}{2\ell^2}\right), \quad \bm{x},\bm{y} \in [-1,1]^d.
\end{equation}
A large value for the length-scale parameter $\ell$ results in smoother Gaussian processes that are more biased towards certain spatial directions. Conversely, a smaller value of $\ell$ results in rougher Gaussian processes that are less spatially biased. We will vary $\ell$ to investigate the effect of the smoothness of the Gaussian process on the result of the low-rank approximation.}


\begin{algorithm}
    \caption{Infinite-dimensional Nystr{\"o}m approximation}
    \label{alg:nystrom}
    \textbf{input:} Non-negative self-adjoint trace-class $\mathcal{A}: L^2(D) \to L^2(D)$, covariance kernel $K: D \times D \to \mathbb{R}$. Target rank $k$, oversampling parameter $p$.\\
    \textbf{output:} Rank $k+p$ Nystr{\"o}m approximation $\widehat{\mathcal{A}}$ to $\mathcal{A}$ in factored form.
    \begin{algorithmic}[1]
        \State Draw a random quasimatrix $\Omega = \begin{bmatrix} \omega_1& \cdots &\omega_{k+p} \end{bmatrix}$ with columns  i.i.d. from $\mathcal{GP}(0,K)$.
        \State Orthonormalize columns of $\Omega$: $Q = \text{orth}(\Omega) = \begin{bmatrix} q_1 & \cdots & q_{k+p}\end{bmatrix}$.\label{line:orthonormalization}
        \State Apply operator $\mathcal{A}$ to $Q$: $Y = \mathcal{A} Q = \begin{bmatrix} \mathcal{A}q_1 & \cdots & \mathcal{A}q_{k+p} \end{bmatrix}$.
        \State $\nu = \epsilon \|Y\|_{\HS}$ where $\epsilon$ is equal to the machine precision.\label{line:nu}
        \State Compute shifted $Y_{\nu} = Y + \nu \Omega$. \label{line:shift}
        \State Compute Cholesky factorization of $\bm{\Omega}^* \bm{Y}_{\nu} = \bm{R}^T \bm{R}$. \Comment{First compute the symmetric part of $\bm{\Omega}^* \bm{Y}_{\nu}$ if needed.} \label{line:chol}
        \State Perform a triangular solve to compute $B = Y_{\nu} \bm{R}^{-1}$.
        \State Compute the Hilbert-Schmidt decomposition of $B = \widehat{U} \bm{S} \bm{V}^*$.
        \State Remove shift $\widehat{\bm{\Sigma}} = \max\{\bm{S}^2 - \nu \bm{I},0\}$, where the maximum is taken entry-wise.
        \State \textbf{return} $\widehat{\mathcal{A}} = \widehat{U}\widehat{\bm{\Sigma}}\widehat{U}^*$ in factored form.
    \end{algorithmic}
\end{algorithm}

\subsection{A pretty function}
\begin{figure}[htbp]
    \centering
    \begin{overpic}[width=\textwidth]{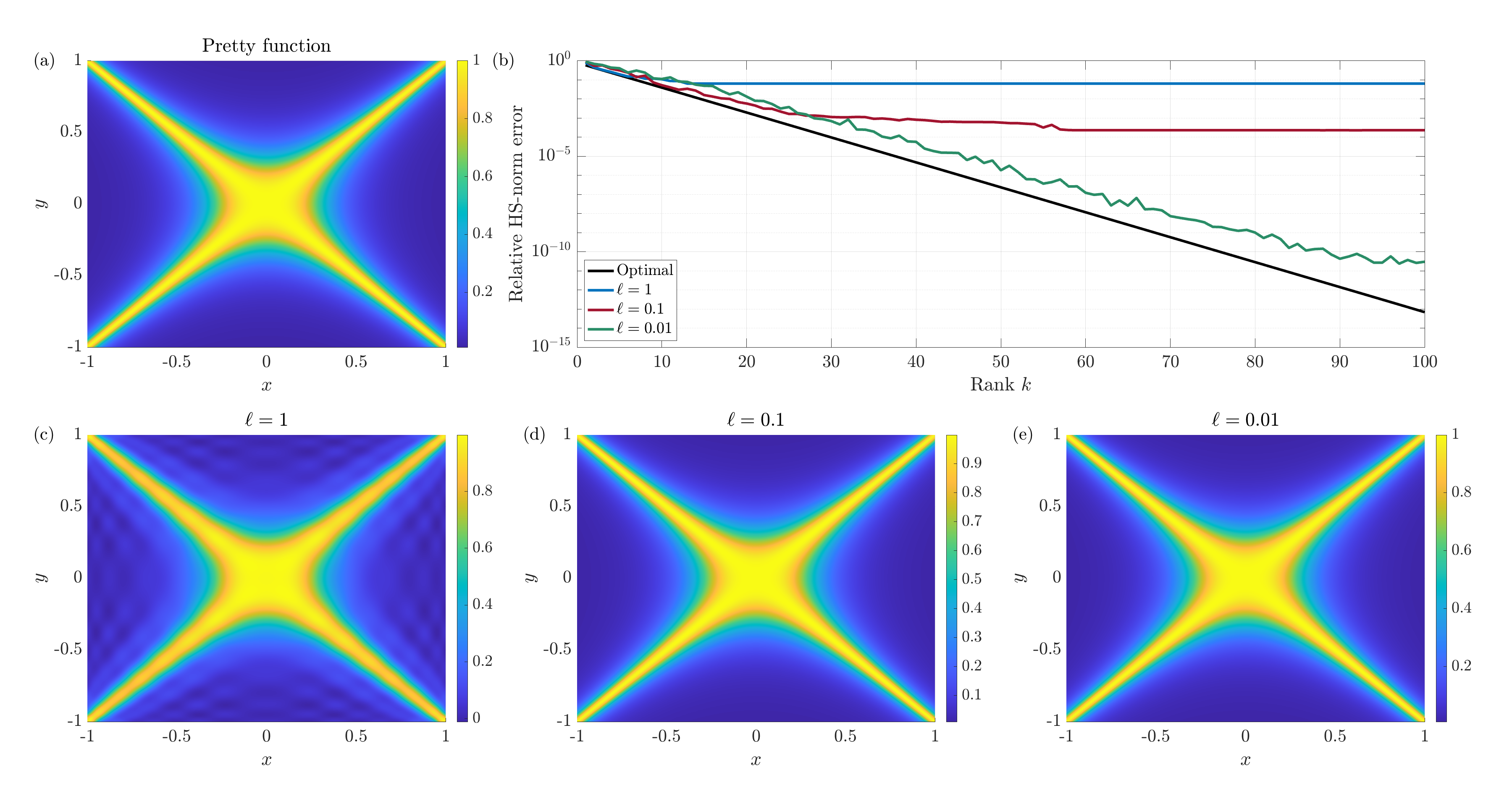}
    \end{overpic}
    \caption{(a) Exact kernel defined by \eqref{eq:pretty} along with convergence of the Nystr{\"o}m approximation for different values of $\ell$ (b). (c)-(e) Rank-$\rev{100}$ Nystr{\"o}m approximations of the kernel for $\ell = 1,0.1,0.01$, respectively.
    }
    \label{fig:pretty}
\end{figure}

In this first example, we compute an approximation to the integral operator \rev{$\mathcal{A}f = \int_{-1}^1 G(\cdot,y) f(y) \d y$} defined by the kernel~\cite{pretty}
\begin{equation}\label{eq:pretty}
    G(x,y) = \frac{1}{1 + 100(x^2 - y^2)^2}, \quad x,y \in [-1,1],
\end{equation}
using \cref{alg:nystrom}. We vary the number of sketches $k$ (i.e., target rank) and report in \cref{fig:pretty} the relative error in the Hilbert-Schmidt norm or, equivalently the $L^2$ norm of $G-G_k$, where $G_k$ is the low-rank Nystr{\"o}m approximant. As can be seen from the figures, setting $\ell = 1$ yields a poor approximation to the kernel where the error stagnates for $k>15$, and $\ell = 0.01$ is better to obtain high-accuracy approximations.

\subsection{Matérn Kernels}
\rev{In this second example, we consider the integral operator \rev{$\mathcal{A}f = \int_{-1}^1 G(\cdot,y) f(y) \d y$} defined by the Matérn-$1/2$, $3/2$, and $5/2$ kernels \cite[Chapt.~4]{gp}:
\begin{align}
     & G_{1/2}(x,y) = \exp(-|x-y|), \quad x,y \in [-1,1];\label{eq:12}                                                                      \\
     & G_{3/2}(x,y) = (1+\sqrt{3}|x-y|) \exp\left(-\sqrt{3}|x-y|\right), \quad x,y \in [-1,1];\label{eq:32}                                 \\
     & G_{5/2}(x,y) = \left(1+\sqrt{5}|x-y| + \frac{5}{3}(x-y)^2\right) \exp\left(-\sqrt{5}|x-y|\right), \quad x,y \in [-1,1].\label{eq:52}
\end{align}
The Matérn class is an important class of low-regularity covariance kernels, frequently appearing in machine learning, where a common task is to compute low-rank approximations of discretizations of the associated kernels. A key benefit of the proposed framework in this work is that once a low-rank approximation to the continuous kernel function is computed, one can obtain a low-rank approximation to \emph{any} discretization of the kernel through a discretization of the low-rank approximation. }

\rev{In this example, we are going to demonstrate how one can incorporate prior knowledge of the kernel into the random fields to obtain more accurate low-rank approximations. }

\begin{figure}[htbp]
    \centering
    \begin{overpic}[width=\textwidth]{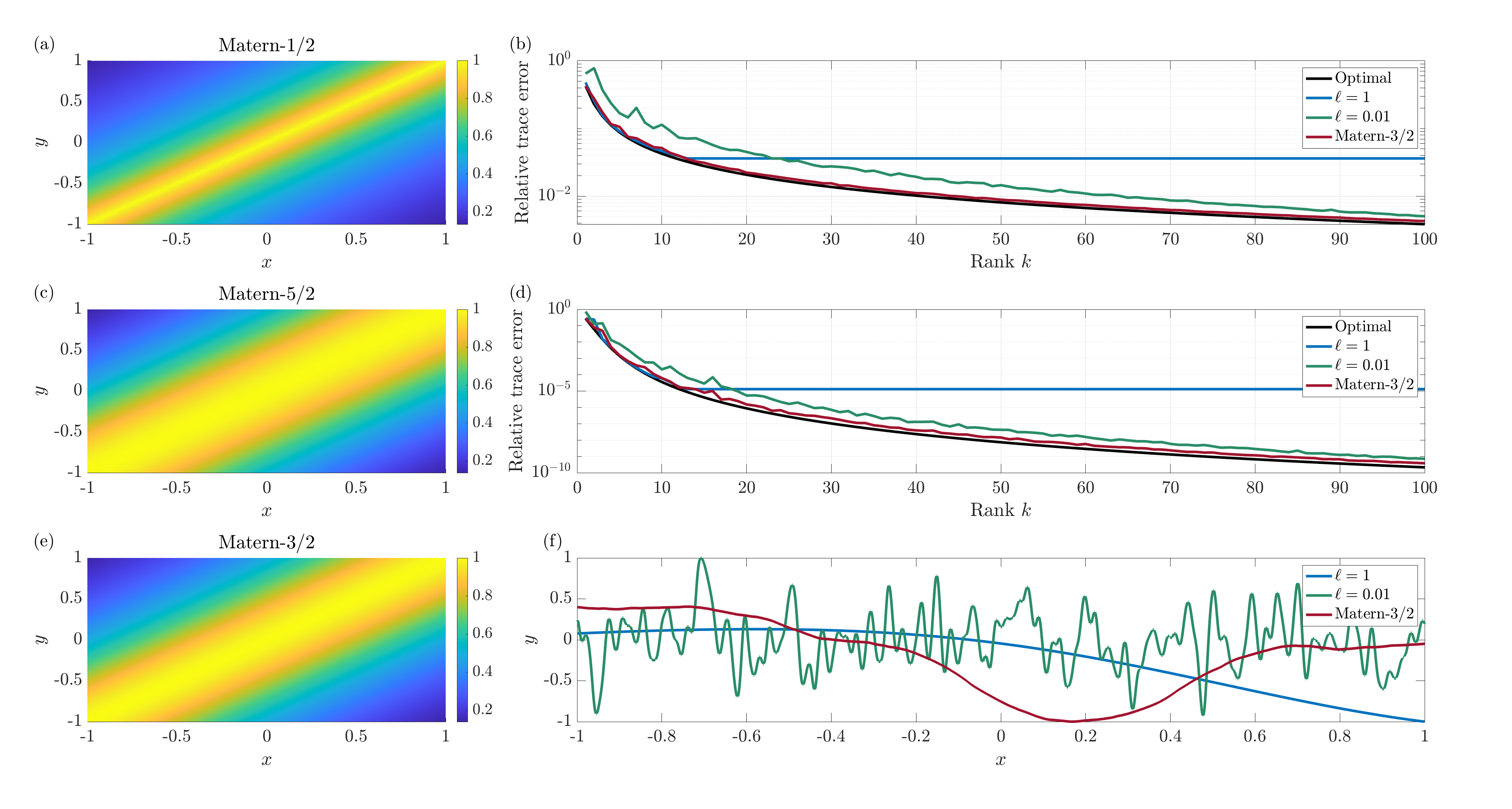}
    \end{overpic}
    \caption{\rev{The contour plots (a),(c), and (e) show the Matérn kernels defined in \eqref{eq:12}-\eqref{eq:52}. The error plots (b) and (d) show the relative error in the trace norm when approximating $G_{1/2}$ and $G_{5/2}$, respectively, when using the squared exponential kernel in \eqref{eq:sekernel} with $\ell = 1, 0.01$ and the Matérn-3/2 kernel in \eqref{eq:32} as covariance kernels for the random fields. \emph{Optimal} denotes the best low-rank approximation error. The plots in (f) show sample paths of the different Gaussian processes. The paths have been normalized so that the maximum absolute value is equal to 1. }}
    \label{fig:matern}
\end{figure}
\begin{figure}[htbp]
    \centering
    \begin{overpic}[width=\textwidth]{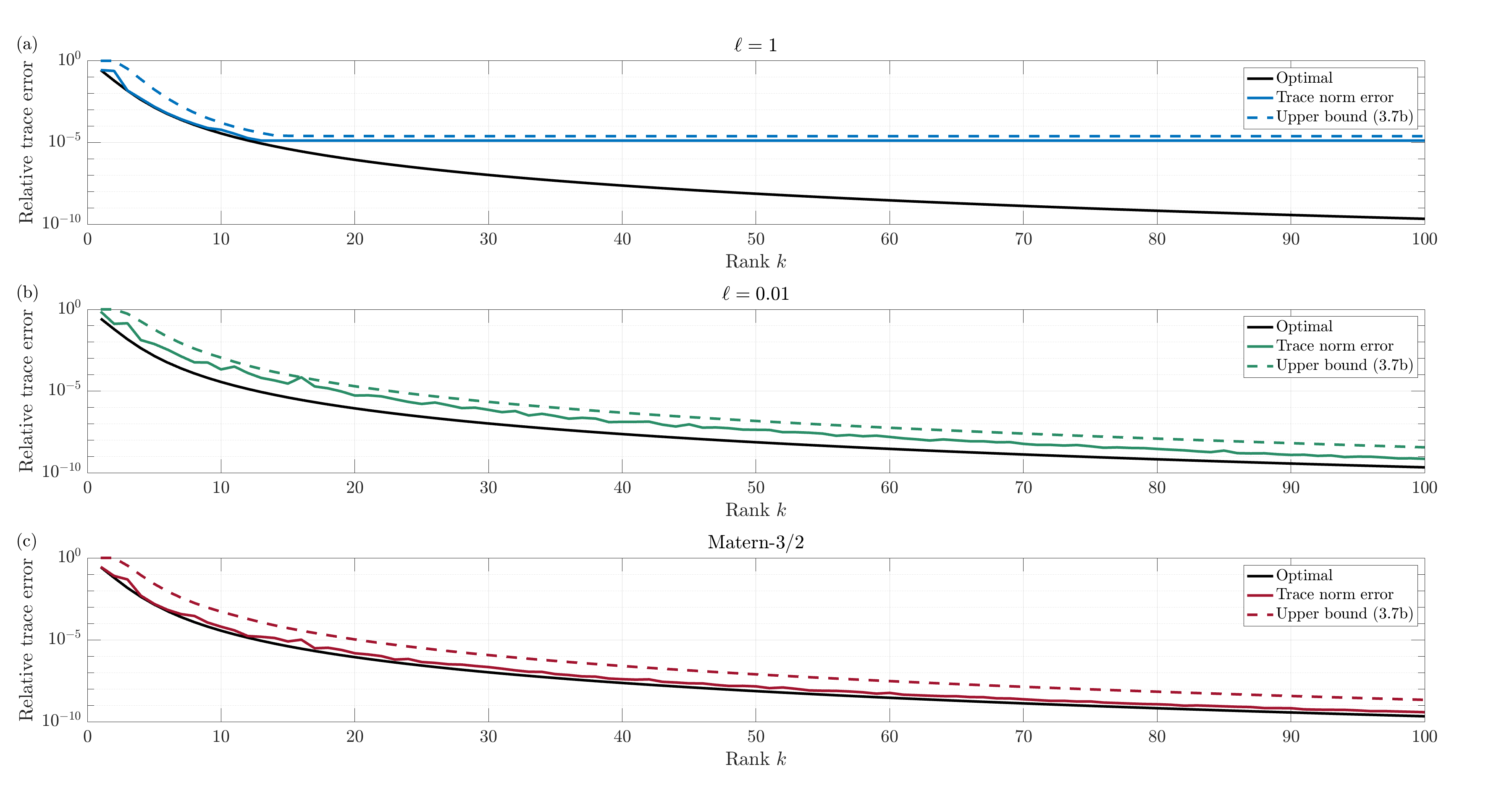}
    \end{overpic}
    \caption{\rev{Figure comparing the numerical results against the theoretical upper bound in \eqref{eq:inf_expectation_tr} when approximating the Matérn-5/2 kernel using the squared-exponential with length parameters $\ell = 1$ and $\ell = 0.01$ and the Matérn-3/2 kernel as covariance kernels. \emph{Optimal} denotes the best low-rank approximation error. The full line shows the numerical results. The dashed lines shows the theoretical upper bound in \eqref{eq:inf_expectation_tr}.}}
    \label{fig:ub}
\end{figure}

\rev{We compute a Nyström approximation for $G_{1/2}$ and $G_{5/2}$ using the squared-exponential kernel with length-scale parameters $\ell = 1$ and $\ell = 0.01$, as well as $G_{3/2}$ as a covariance kernel for the random fields. Using $G_{3/2}$ as a covariance kernel incorporates prior knowledge of the operator's eigenfunctions, since the Matérn-$\nu$ kernels vary continuously with $\nu$. The results are presented in \Cref{fig:matern}, which demonstrate that incorporating such prior knowledge can lead to more accurate low-rank approximations. For the Matérn-5/2 kernel, we have also included a figure comparing the numerical results with the upper bound in \eqref{eq:inf_expectation_tr}; see \Cref{fig:ub}. Notably, from \Cref{fig:matern} and \Cref{fig:ub}, we observe that using $G_{3/2}$ as a covariance kernel for random fields yields excellent results. For the squared-exponential kernel, when $\ell$ is small, there is significantly less variance in the directions of $U_2$, i.e., \(\beta_k^{(\Tr)}\) is small for the corresponding random field. This results in near-optimal low-rank approximations for small values of $k$. However, as $k$ increases, the random fields become too spatially biased, leading to large values of $\delta_k^{(k)}$. In contrast, when $\ell$ is small, the random fields exhibit a higher variance in all directions due to less spatial bias. This leads to larger values of $\beta_k^{(\Tr)}$ for small $k$, but as higher-frequency terms are needed to capture the eigenfunctions associated with smaller eigenvalues, the approximation becomes more efficient for larger $k$. } 

\rev{In cases where no prior knowledge of the operator is available, we recommend using the roughest possible random fields, subject to computational constraints. Although sampling from rougher random fields is typically more computationally expensive than sampling from smoother ones, since the corresponding Karhunen-Loève expansion converges slower, they offer a safer choice, as they are less spatially biased. }



\subsection{Approximating a Karhunen-Loève expansion}
\rev{In our next example we consider the task of sampling from a Gaussian process. Given a symmetric positive semi-definite kernel $G: D \times D \to \mathbb{R}$, inducing an integral operator $\mathcal{A}f = \int_D G(\cdot,y) f(y) \d y$, a common task in computational statistics is to sample from $\mathcal{GP}(0,G)$. This can be achieved in a few ways. For example, if the Mercer expansion of $G$ is known, one can sample from $\mathcal{GP}(0,G)$ through the Karhunen-Loève expansion \eqref{eq:klexpansion}:
\begin{equation*}
    \omega = \sum\limits_{i=1}^{\infty} \lambda_i^{1/2} \zeta_iu_i \sim \mathcal{GP}(0,G),
\end{equation*}
where $\zeta_1,\zeta_2,\cdots \sim \mathcal{N}(0,1)$ are mutually independent. However, computing a Mercer expansion of $G$ can be computationally intractable. One can also sample from a discretized Gaussian process by sampling from $\mathcal{N}(\bm{0}, \bm{G})$, where $\bm{G}_{ij} = G(x_i,x_j)$ and $\{x_i\}_{i=1}^n$ are the discretization points. This requires computing a square-root or a Cholesky factorization of $\bm{G}$, which costs $O(n^3)$ operations and becomes prohibitively expensive for large $n$.}

\rev{If we have access to an accurate Nyström approximation $G_k$ to $G$ one can easily sample from $\mathcal{GP}(0,G_k)$ as
\begin{equation*}
    \widehat{\omega} = \sum\limits_{i=1}^{k} \widehat{\lambda}_i^{1/2} \widehat{\zeta}_i\widehat{u}_i\sim \mathcal{GP}(0,G_k),
\end{equation*}
where $(\widehat{\lambda}_i,\widehat{u}_i)$ are the eigenpairs of $G_k$ for $i = 1,\ldots,k$, and $\widehat{\zeta}_1,\ldots,\widehat{\zeta}_k \sim \mathcal{N}(0,1)$ are independent. Computing a Nyström approximation of $\mathcal{G}$ is usually significantly less expensive than computing a Mercer expansion of $G$. Furthermore, choosing $\widehat{\zeta}_i$ so that $\mathbb{E}[\zeta_i \widehat{\zeta}_j] = \langle u_i,\widehat{u}_j \rangle_{L^2(D)}$ yields $\mathbb{E}\|\omega - \widehat{\omega}\|_{L^2(D)}^2 \leq \|G-G_k\|_{\Tr}$.}\footnote{\rev{This fact requires that $G-G_k$ is a positive semi-definite kernel, which is true whenever $G_k$ is the kernel of the Nyström approximation. However, this fact is not true for arbitrary approximations.}} \rev{Hence, $\mathcal{W}_2^2(\mathcal{GP}(0,G),\mathcal{GP}(0,G_k)) \leq \|G-G_k\|_{\Tr}$, where $\mathcal{W}_2$ denotes the Wasserstein distance. Therefore, as long as $G_k$ is an accurate low-rank approximation to $G$, $\mathcal{GP}(0,G_k)$ is a good approximation to $\mathcal{GP}(0,G)$ in the Wasserstein distance. }

\begin{figure}[htbp]
    \centering
    \begin{overpic}[width=\textwidth]{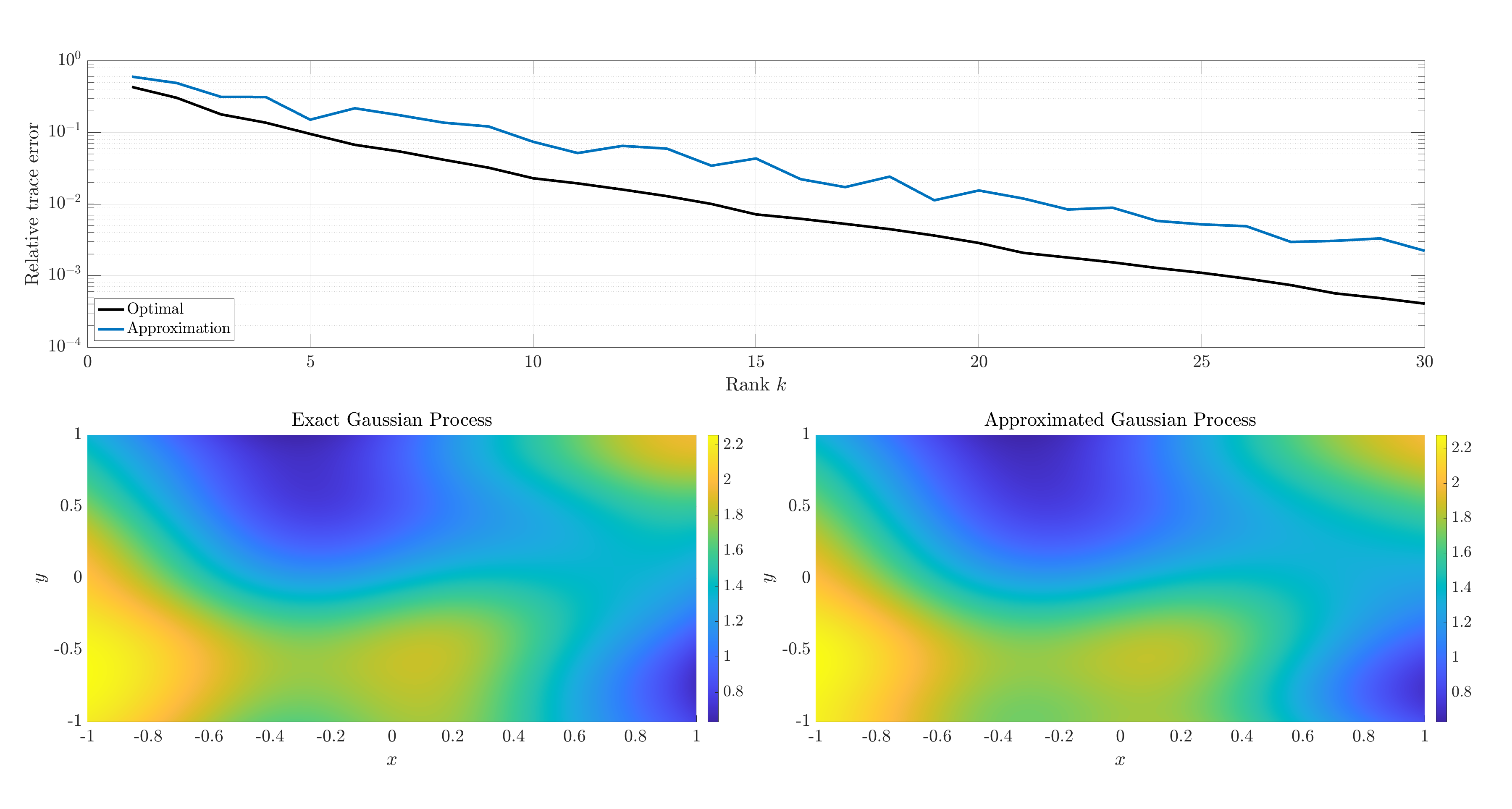}
    \end{overpic}
    \caption{\rev{The left contour plot shows an exact sample from $\mathcal{GP}(0,G)$ and the left contour plot shows a sample from $\mathcal{GP}(0,G_{30})$ so that the mean-squared error equals $\|G-G_{30}\|_{\Tr}$. The error plots show the relative error in the trace norm. \emph{Optimal} denotes the best low-rank approximation error.}}
    \label{fig:gp}
\end{figure}

\rev{In this example we consider $G$ to be the squared-exponential kernel \eqref{eq:sekernel} with $d = 2$ and $\ell = 0.4$. The random fields are drawn from $\mathcal{GP}(0,K)$ with $K(\bm{x},\bm{y}) = \sum_{i=1}^{625} p_i(\bm{x}) p_i(\bm{x})$ where $\{p_1,\ldots,p_{625}\}$ is an orthonormal basis obtained by taking the tensor product of the first 25 Legendre polynomials. For the purpose of this experiment, we discretize the squared-exponential kernel by expanding it onto the orthonormal basis obtained by taking the tensor product of the first $100$ Legendre polynomials. The results are displayed in \cref{fig:gp}. 
}

\subsection{An example from Bayesian inverse problems}
\rev{Our final example comes from A-optimal design of Bayesian inverse problems. In this setting, we have a linear parameter-to-observable map $\mathcal{F}: L^2(D) \to \mathbb{R}^n$, which maps, for example, the initial condition of a partial differential equation to the function values at a set of discrete points of the solution. In the context of A-optimal design of infinite-dimensional linear inverse problems one wants to compute a low-rank approximation of the operator
\begin{equation*}
    \mathcal{A} = \mathcal{C}_0^{1/2} \mathcal{F}^* \bm{\Sigma}_{\text{noise}} \mathcal{F}\mathcal{C}_0^{1/2},
\end{equation*}
where $\mathcal{C}_0$ is the covariance operator of the Gaussian prior distribution and $\bm{\Sigma}_{\text{noise}}$ is the covariance matrix of measurement errors \cite{alexanderian2014optimal,saibabarandom}.}

\rev{As performed in \cite{alexanderian2014optimal,saibabarandom}, we let $\mathcal{C}_0$ be the solution operator of the Biharmonic equation on $\Omega=[-1,1]^2$ with Dirichlet boundary conditions, $\bm{\Sigma}_{\text{noise}} \propto \bm{I}$, and denote by $\mathcal{F}$ the operator that maps the (space-dependent) initial condition $\theta \in L^2(\Omega)$ to the time-dependent solution $u$ of the following parabolic differential equation:
\begin{align}
\begin{split}
    \partial_t u = \Delta u + 0.1 u \quad \text{in } \Omega \times \mathbb{R}_{\geq 0},\\
    u = 0 \quad\text{on }\partial \Omega,\quad\text{and}\quad
    u(\cdot,0) = \theta\quad\text{in } \Omega,
    \end{split}\label{eq:pde}
\end{align}
on a $50 \times 50$ equispaced grid of $\Omega$ and at times $t = 0.1, 0.5,1$. }

\begin{figure}[htbp]
    \centering
    \begin{overpic}[width=\textwidth]{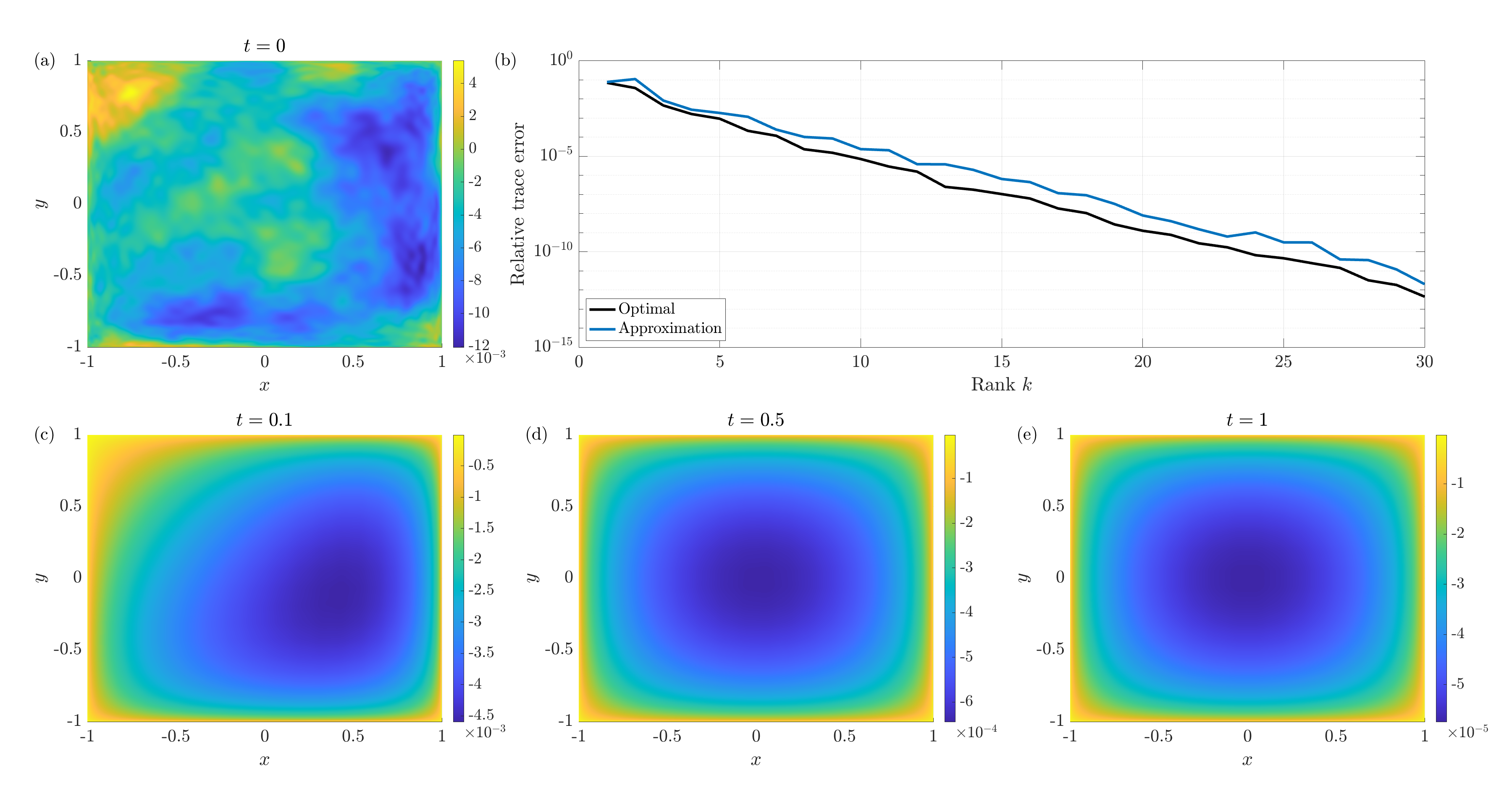}
    \end{overpic}
    \caption{\rev{The contour plots shows the solution of \eqref{eq:pde} for different times. The error plots show the relative error in the trace norm norm. \emph{Optimal} denotes the best low-rank approximation error.}}
    \label{fig:inverse}
\end{figure}

\rev{For the purpose of this experiment we solve \eqref{eq:pde} by first discretizing in space using finite differences on a $50 \times 50$ equispaced grid, yielding an ordinary differential equation (ODE). We then proceed in time by solving the resulting ODE using exponential integrators. We let the covariance kernel for the random fields be the solution operator of the Biharmonic equation on $[-1,1]^2$ with Direchlet boundary conditions. The results are displayed in \cref{fig:inverse}, which demonstrates that the Nyström approximation yields near-optimal low-rank approximations to the operator. }

\section{Conclusions}
\label{section:conclusions}

In this work, we have presented an infinite-dimensional extension of the randomized Nystr{\"o}m approximation. We first established bounds for the finite-dimensional Nystr{\"o}m approximation when the columns of the sketch matrix $\bm{\Omega}$ are drawn independently from a non-standard Gaussian distribution $\mathcal{N}(\bm{0},\bm{K})$. Subsequently, through a continuity argument, we provided analogous bounds for the infinite-dimensional extension of the Nystr{\"o}m approximation for Hilbert-Schmidt operators. Additionally, in the process of analyzing the Nystr{\"o}m approximation for Hilbert-Schmidt operators, we have also improved the existing bounds for the randomized SVD for Hilbert-Schmidt operators. \rev{Our analysis is complemented by a number of numerical experiments on the application of the proposed infinite-dimensional Nystr\"om method and demonstrate its potential in Gaussian process sampling or Bayesian inverse problems.}

\section*{Data availability}

All experiments have been performed in MATLAB (version 2020a) on a MacBook Pro with a 2.3 GHz Intel Core i7 processor with 4 cores, and the scripts to reproduce all figures in this paper are publicly available at \url{https://github.com/davpersson/infNystrom}.

\appendix

\section{Properties of Gaussian matrices}\label{sec_prop_Gaussian} The following lemma is a consequence of the symmetry of standard normal random variables.

The following lemma summarizes results on the expected norms of scaled and shifted Gaussian matrices for the Frobenius, Schatten-4, and spectral norms.

\begin{lemma}[Expected norm of shifted Gaussian matrices] \label{lemma:shifted_4norm}
    Let $\bm{C} \in \mathbb{R}^{m_1 \times m_2}, \bm{D} \in \mathbb{R}^{n_1 \times n_2}$, and $\bm{B} \in \mathbb{R}^{m_1 \times n_2}$ be three arbitrary matrices, and consider an $m_2 \times n_1$ standard Gaussian matrix $\bm{\Psi}$. Then, the following relations hold:
    \begin{subequations}
        \begin{align}
            \mathbb{E}\big[\|\bm{B} + \bm{C}\bm{\bm{\Psi}}\bm{D}\|_\Frob^2\big] & = \|\bm{B}\|_{\Frob}^2 + \|\bm{C}\|_\Frob^2 \|\bm{D}\|_\Frob^2, \label{eq_exp_shift_Frob}                                                           \\
            \mathbb{E}\big[\|\bm{C}\bm{\Psi}\bm{D}\|_{(4)}^4\big]               & = \|\bm{C}\|_{(4)}^4\|\bm{D}\|_{(4)}^4 + \|\bm{C}\|_{\Frob}^4\|\bm{D}\|_{(4)}^4 + \|\bm{C}\|_{(4)}^4\|\bm{D}\|_{\Frob}^4, \label{eq_exp_shift_four} \\
            \mathbb{E}\big[\|\bm{C}\bm{\bm{\Psi}}\bm{D}\|_2^2\big]              & \leq \left(\|\bm{C}\|_\Frob\|\bm{D}\|_2+\|\bm{C}\|_2\|\bm{D}\|_\Frob\right)^2, \label{eq_exp_shift_two}                                             \\
            \mathbb{E}\big[\|\bm{C} \bm{\Psi} \bm{D}\|_\Frob^4]                 & = 2\|\bm{C}\|_{(4)}^4 \|\bm{D}\|_{(4)}^4 + \|\bm{C}\|_{\Frob}^4 \|\bm{D}\|_{\Frob}^4.\label{eq:exp_4power}
        \end{align}
    \end{subequations}
\end{lemma}
\begin{proof}
    We introduce the matrix $\bm\Phi = \bm C \bm \Psi \bm D$ and begin by proving \cref{eq_exp_shift_Frob}. Using the linearity of trace and expectation we have:
    \[\mathbb{E}[\|\bm B + \bm\Phi\|_\Frob^2]=\|\bm B\|_\Frob^2 + \mathbb{E}[\|\bm\Phi\|_{\Frob}^2]+2\mathbb{E}[\tr(\bm B^* \bm\Phi)]=\|\bm{B}\|_{\Frob}^2 + \|\bm{C}\|_\Frob^2 \|\bm{D}\|_\Frob^2,\]
    where we combined~\cite[Prop.~10.1]{rsvd} with the equality $\mathbb{E}[\tr(\bm B^* \bm\Phi)]=0$ since $\bm \Psi$ has zero mean. A similar result is found in \cite[Lem.~3.11]{diouane2022general}.
    The equality~\cref{eq_exp_shift_four} and inequality~\cref{eq_exp_shift_two} can be found in~\cite[Lem.~B.1]{tropp2023randomized} and~\cite[Prop.~B.3]{frangella2021randomized}, respectively.
    We now conclude with the proof of \eqref{eq:exp_4power}. Let $\bm{E} = \bm{D}^* \otimes \bm{C}$ and $\bm{\psi} = \text{vec}(\bm{\Psi})$. Then,
    \begin{align*}
        \mathbb{E}[\|\bm{C}\bm{\Psi} \bm{D}\|_\Frob^4] & = \mathbb{E}[(\bm{\psi}^* \bm{E}^*\bm{E} \bm{\psi})^2] = \Var(\bm{\psi}^* \bm{E}^* \bm{E} \bm{\psi}) + (\mathbb{E}[\bm{\psi}^* \bm{E}^* \bm{E} \bm{\psi}])^2 \\ & = 2\|\bm{E}^* \bm{E}\|_\Frob^2 + \tr(\bm{E}^* \bm{E})^2 = 2\|\bm{C
        }\|_{(4)}^4 \|\bm{D}\|_{(4)}^4 + \|\bm{C}\|_{\Frob}^4 \|\bm{D}\|_{\Frob}^4,
    \end{align*}
    where we used a standard result that the variance of a quadratic form with Gaussian random vectors is $2\|\bm{E}^T \bm{E}\|_\Frob^2$,
\end{proof}


\begin{lemma} \label{lemma:inversewishart_norm_expectation}
    Let $\bm\Omega_1\in \mathbb{R}^{k\times(k+p)}$ be a random matrix whose columns are i.i.d. $\mathcal{N}(0,\tK_{11})$ and let $\bm B\in \R^{k\times n}$ be an arbitrary matrix. Then, the following relation holds for $k \geq 1$ and $p\geq 2$:
    \begin{subequations}
        \begin{align}
            \mathbb{E}[\|\bOmega_1^{\dagger}\bm{B}\|_\Frob^2] & = \frac{\tr(\bm{B}^* \tK_{11}^{-1} \bm{B})}{p-1} = \frac{\|\bm{B}^* \tK_{11}^{-1/2}\|_\Frob^2}{p-1}, \label{eq_inv_wish_1} 
        \end{align}
    \end{subequations}
    Additionally, for $p, k \geq 2$ we have:
    \begin{subequations}
        \begin{align}
            \mathbb{E}[\|\bOmega_1^{\dagger}\|_2^{2}] & \leq \frac{e^{2}(k+p)}{(p-1)(p+1)}\|\tK_{11}^{-1}\|_2. \label{eq_inv_wish_2}
        \end{align}
    \end{subequations}
    Finally, for $k \geq 1$ and $p\geq 4$ we have
    \begin{subequations}
        \begin{align}
            \mathbb{E}[\|\bOmega_1^{\dagger}\|_{(4)}^4] & = \frac{(p-1)\|\tK_{11}^{-1}\|_\Frob^2 + \tr(\tK_{11}^{-1})^2}{p(p-1)(p-3)}\label{eq_inv_wish_3} \le k \frac{k+p-1}{p(p-1)(p-3)} \|\tK_{11}^{-1}\|_2^2,      \\
            \mathbb{E}[\|\bOmega_1^{\dagger}\|_\Frob^4] & = \frac{(p-2)\tr(\tK_{11}^{-1})^2 + 2 \|\tK_{11}^{-1}\|_\Frob^2}{p(p-1)(p-3)} \le k \frac{kp-2k+2}{p(p-1)(p-3)} \|\tK_{11}^{-1}\|_2^2. \label{eq_inv_wish_4}
        \end{align}
    \end{subequations}
\end{lemma}
\begin{proof}
    First, \cref{eq_inv_wish_1} is proven in~\cite[Lem.~3.12]{diouane2022general} and \cref{eq_inv_wish_2} follows from in~\cite[Lem.~3.1]{nakatsukasa2020fast} and the fact $\mathbb{E}[\|\bm{\Omega}_1^{\dagger}\|_2^2] \leq \mathbb{E}[\|\bm{\Psi}^{\dagger}\|_2^2]\|\bm{K}_{11}^{-1}\|_2$ where $\bm{\Psi}$ is a $k \times (k+p)$ standard Gaussian matrix.
    We proceed to prove the equality in~\cref{eq_inv_wish_3} and introduce the random variable $\bm{X} = \bOmega_1\bOmega_1^* \sim \mathcal{W}_k(\tK_{11},k+p)$, such that
    \[\mathbb{E}\big[\|\bOmega_1^{\dagger}\|_{(4)}^4\big] = \mathbb{E}\big[\|\bm{X}^{-1}\|_\Frob^2\big] = \mathbb{E}\big[\tr(\bm{X}^{-2})\big] = \frac{\|\tK_{11}^{-1}\|_\Frob^2}{p(p-3)} + \frac{\tr(\tK_{11}^{-1})^2}{p(p-1)(p-3)},\]
    where the last equality follows from~\cite[Thm.~2.4.14]{kollorosen}. The equality in~\cref{eq_inv_wish_4} follows from the fact that $\mathbb{E}[\|\bOmega_1^{\dagger}\|_\Frob^4] = \mathbb{E}[\tr(\bm{X}^{-1})^2] = \mathbb{E}[\tr\left(\bm{X}^{-1} \otimes \bm{X}^{-1}\right)]$. We then exploit a relation on $\mathbb{E}[\bm{X}^{-1} \otimes \bm{X}^{-1}]$~\cite[Thm.~2.4.14]{kollorosen} to obtain
    \begin{align*}
        \mathbb{E}\big[\|\bOmega_1^{\dagger}\|_\Frob^4\big] & = \frac{(p-2)\tr(\tK_{11}^{-1} \otimes \tK_{11}^{-1}) + \tr(\vect(\tK_{11}^{-1})\vect(\tK_{11}^{-1})^*) + \tr(\bm{C}_{k \times k} (\tK_{11}^{-1} \otimes \tK_{11}^{-1})) }{p(p-1)(p-3)} \\
                                                            & = \frac{(p-2)\tr(\tK_{11}^{-1})^2 + 2 \|\tK_{11}^{-1}\|_\Frob^2}{p(p-1)(p-3)},
    \end{align*}
    where the second equality comes from the relation $\tr(\bm{C}_{k \times k} (\tK_{11}^{-1} \otimes \tK_{11}^{-1})) = \|\tK_{11}^{-1}\|_\Frob^2$, where $\bm{C}_{k \times k}$ is the commutation matrix.\footnote{The commutation matrix $\bm{C}_{k \times k} \in \mathbb{R}^{k^2 \times k^2}$ is a $k \times k$ block matrix, with blocks of size $k \times k$. The $(i,j)$-block of $\bm{C}_{k \times k}$ is the matrix $\bm{E}_{ji}$ with entries $(\bm{E}_{ji})_{k\ell} = \delta_{jk}\delta_{i\ell}$.} Finally, the inequalities in~\cref{eq_inv_wish_3} and~\cref{eq_inv_wish_4} follow from standard norm inequalities.
\end{proof}
The next lemma is a generalization and a consequence of~\cite[Prop.~10.4]{rsvd}, which provides tail bounds on the Frobenius and spectral norms of pseudoinverted standard Gaussian matrices.

\begin{lemma}[Norm bounds for a pseudoinverted scaled Gaussian matrix]\label{lemma:omega1_tailbound}
    Let $\bm\Omega_1\in \mathbb{R}^{k\times(k+p)}$ be a Gaussian matrix whose columns are i.i.d.~and follow a multivariate Gaussian distribution with mean zero and covariance matrix $\tK_{11}$. Then, the following relations hold for $p\geq 4$ and all $t\geq 1$:
    \begin{subequations}
        \begin{align}
            \mathbb{P}\left\{\|\bOmega_1^{\dagger}\|_\Frob > \sqrt{\frac{3\tr(\tK_{11}^{-1})}{p+1}} t\right\}             & \leq t^{-p}, \label{eq_omega_1_frob}         \\
            \mathbb{P}\left\{\|\bOmega_1^{\dagger}\|_{2} > \frac{e \sqrt{(k+p) \|\tK_{11}^{-1}\|_2}}{p+1} t\right\}       & \leq t^{-(p+1)}, \label{eq_omega_1_spec}     \\
            \mathbb{P}\left\{\|\bOmega_1^{\dagger}\|_{(4)} > \frac{e \sqrt{(k+p) \|\tK_{11}^{-1}\|_\Frob}}{p+1} t\right\} & \leq t^{-(p+1)}. \label{eq_omega_1_schatten}
        \end{align}
    \end{subequations}
\end{lemma}
\begin{proof}
    Note that~\cref{eq_omega_1_frob} is a restatement of~\cite[Lem.~3]{boulle2022learning}. Because $\bOmega_1 = \tK_{11}^{1/2}\bm\Psi$, where $\bm\Psi$ is a standard Gaussian matrix, it follows that
    $\|\bm\Psi^\dagger \tK_{11}^{-1/2}\|_{(s)}\leq \|\bm\Psi^\dagger\|_2\|\tK_{11}^{-1/2}\|_{(s)}$ for any Schatten-$s$ norm. Moreover, the combination with \cite[Prop.~10.4]{rsvd} implies
    \[
        \mathbb{P}\left\{\|\bm\Psi^{\dagger}\|_{2} > \frac{e \sqrt{(k+p) }}{p+1} t\right\}\leq t^{-(p+1)},
    \]
    which yield the bounds~\cref{eq_omega_1_spec,eq_omega_1_schatten} using
    $\|\tK_{11}^{-1/2}\|^2_{2} = \|\tK_{11}^{-1}\|_2$ and $\|\tK_{11}^{-1/2}\|^2_{(4)} = \|\tK_{11}^{-1}\|_\Frob$,
    respectively.
    %
\end{proof}


\bibliographystyle{siamplain}
\bibliography{bibliography}

\end{document}


\maketitle

\section{A detailed example}

Here we include some equations and theorem-like environments to show
how these are labeled in a supplement and can be referenced from the
main text.
Consider the following equation:
\begin{equation}
  \label{eq:suppa}
  a^2 + b^2 = c^2.
\end{equation}
You can also reference equations such as \cref{eq:matrices,eq:bb} 
from the main article in this supplement.

\begin{theorem}
An example theorem.
\end{theorem}
 
\begin{lemma}
An example lemma.
\end{lemma}

Here is an example citation:

\section[Proof of Thm]{Proof of \cref{thm:bigthm}}
\label{sec:proof}

\section{Additional experimental results}
\Cref{tab:foo} shows additional
supporting evidence. 

\begin{table}[htbp]
\footnotesize
  \caption{Example table.}  \label{tab:smfoo}
\begin{center}
  \begin{tabular}{|c|c|c|} \hline
   Species & \bf Mean & \bf Std.~Dev. \\ \hline
    1 & 3.4 & 1.2 \\
    2 & 5.4 & 0.6 \\ \hline
  \end{tabular}
\end{center}
\end{table}

\bibliographystyle{siamplain}
\bibliography{bibliography}